\documentclass{amsart}

\usepackage{graphicx}

\def\R{{\mathbb R}}

\def\N{{\mathbb N}}

\def\Z{{\mathbb Z}}
\def\1{{1\!\!\!1}}
\def\a{{\alpha}}
\def\E{{\mathbb E}}
\def\P{{\mathbb P}}

\def\cal{\mathcal}

\def\ol{\overline}

\def\supp{{\rm{supp}}}
\def\eps{\varepsilon}

\newcommand{\be}{\begin{equation}}
\newcommand{\ee}{\end{equation}}
\numberwithin{equation}{section}

\newtheorem{theorem}{Theorem}
\newtheorem{prop}{Proposition}[section]
\newtheorem{cor}{Corollary}[section]
\newtheorem{defi}{Definition}[section]
\newtheorem{lemma}{Lemma}[section]

\title{Martin boundary of a killed random walk on a half-space}
\author{Irina Ignatiouk-Robert}
\address{
{Universit\'e de Cergy-Pontoise,}
{D\'epartement de math\'ematiques,}
{2, Avenue Adolphe Chauvin,}
{95302 Cergy-Pontoise Cedex,}
{France}}
\date{\today}
\email{Irina.Ignatiouk@math.u-cergy.fr}

\begin{document}
\begin{abstract}  A complete representation of the Martin boundary of killed random walks  on a half-space
$\Z^{d-1}\times\N^*$ is obtained.  In particular, it is proved that the corresponding Martin
  boundary is homemorphic to the half-sphere ${\cal S}^d_+ = \{z\in\R^{d-1}\times\R_+ :
  |z|=1\}$. The method  is based on a combination of ratio
  limits theorems and large deviation techniques. 
\end{abstract}
\maketitle

\section{Introduction}\label{sec1}

The concept of Martin compactification was first introduced for Brownian motion by
Martin~\cite{Martin}. For countable Markov chains with discrete time, the abstract construction
of the Martin compactification was given by Doob~\cite{Doob} and Hunt~\cite{Hunt}, see
also Dynkin~\cite{Dynkin:01}  and 
Rogers and Williams~\cite{Rogers:05}. The main results of this
theory are the following~:

For a transient Markov chain $(Z(n))$ on a countable set $E$ having Green's function
$G(z,z')$, the 
Martin kernel $K(z,z')$ is defined by 
\[
K(z,z') ~=~ G(z,z')/G(z_0,z') 
\]
where $z_0$ is a given reference point in $E$. A sequence $z_n\in E$ is said to converge to a
point on a Martin  boundary $\partial E_M$ of $E$ if it leaves  every finite subset on $E$
and  the sequence  of functions  $K(\cdot,z_n)$  converges point-wise.  According to  this
definition, the Martin  compactification $E_M$ is the unique  smallest compactification of
the  set $E$ for  which the  Martin kernels $K(z,\cdot)$  extend continuously. 

The minimal
Martin boundary $\partial_m E_M$ is the set of all those
$\gamma\in\partial E_M$ for which the function $K(\cdot, \gamma)$ is 
minimal harmonic. Recall that a function $h : E\to\R_+$ is harmonic for $(Z(n))$
if $
\E_z(h(Z(1))) ~=~ h(z)$ 
for all $z\in E$. A harmonic function $h : E\to\R_+$ is minimal if the inequality 
$
0 \leq h' \leq h$ for any other harmonic function $h'$ implies that $h' = ch$ with some
$c>0$. 

By the Poisson-Martin
representation  theorem, for  every  non-negative  harmonic function  $h$  there exists  a
positive Borel measure $\nu$ on $\partial_m E_M$ such that
\[
h(z) = \int_{\partial_m E_M} K(z,\eta) \, d\nu(\eta) 
\]
By Convergence theorem, the sequence
$(Z(n))$ converges $\P_z$ almost surely for every $z\in E$ to a $\partial_m E_M$ valued
random variable. The Martin boundary provides therefore the non-negative harmonic functions and 
shows how the Markov chain $(Z(n))$ goes to infinity.  A good introduction to the theory of
Martin compactification for countable Markov chains and a wide literature of
related results   is given in the book of Woess~\cite{Woess}.

An explicit description  of the Martin compactification is  usually a non-trivial problem.
A  large number  of  results  in this  domain  has been  obtained  for homogeneous  random
walks. For homogeneous  random walks on $\Z^d$, the Martin boundary  was identified by Ney
and  Spitzer~\cite{Ney-Spitzer}. They  considered an  irreducible   random walk
$(Z(t))$ on $\Z^d$ with transition probabilities
$p(z,z') ~=~ \mu(z'-z)$ 
and a non zero mean,   for which the jump generating function 
\be\label{e1-0}
\varphi(a) ~\dot=~ \sum_{z\in \Z^d} \mu(z)
~\exp(a\cdot z)
\ee 
is finite in a neighborhood of the set 
\[
D ~\dot=~ \{a\in \R^d :
\varphi(a) \leq 1\}.
\]
 For
such a random walk, the set $D$ is compact and convex, the gradient  $
\nabla \varphi(a)$
exists everywhere on $\R^d$ and does not vanish on the boundary $\partial D \dot= \{a :
\varphi(a) = 1\}$, and the mapping
\be\label{e1-1}
 q (a) = {\nabla \varphi(a)}/{|\nabla \varphi(a)|}
\ee
determines a homeomorphism between  $\partial D$ and the unit sphere ${\cal S}^d$ in $\R^d$
(see~\cite{Hennequin}). Using  exponential change of measure and the local limit
theorem, Ney and Spitzer calculated the exact asymptotics of the Green's function  
\[
G(z,z') = \sum_n \P_z(Z(n)=z')
\]  
and deduced that 
for any $ a\in\partial D$,  and any sequence of points $z_n\in\Z^d$ with
$\lim_{n} |z_n|\to\infty$ and $\lim_{n} z_n/|z_n| = q(a),$ 
\[
{G(z,z_n)}/{G(0,z_n)} ~\to~ \exp(a\cdot z) \;\text{ as \; $n\to\infty$, } \quad \quad \forall \; z\in\Z^d. 
\]
Hence, for the  homogeneous random walk on $\Z^d$, a sequence of
points  $z_n\in\Z^d$ with  $\lim_n|z_n|  = \infty$  converges  to a  point  on the  Martin
boundary if and only if the sequence $z_n/|z_n|$  converges to a point on a unit sphere in
$\R^d$. The  Martin compactification of the  lattice $\Z^d$ determined  by the homogeneous
random  walk  is  therefore  homeomorphic  to   the  closure  of  the  set  $\left\{  w  =
{z}/{(1+|z|)}~:~ z\in\Z^d\right\}$ in $\R^d$.

For a wide literature of results where the Martin boundary was identifies
for more general homogeneous Markov chains such as random walks on free groups, hyperbolic
graphs, Cartesian products we refer to  the  book of
Woess~\cite{Woess} and the references therein. 

\medskip 

Only few results identify the Martin boundary for non-homogeneous Markov
 chains. For random walks on
 non-homogeneous trees the Martin boundary was described by Cartier~\cite{Cartier}.  
 Doney~\cite{Doney:02} identified the harmonic functions and the Martin boundary of 
a homogeneous random walk $(Z(n))$ on $\Z$ killed  on the negative half-line
$\{z : z<0\}$. For space-time random  walk $S(n)
=(Z(n),n)$ for a  homogeneous random walk $Z(n)$ on $\Z$ killed  on the negative half-line
$\{z : z<0\}$ the Martin boundary was described by  Alili and
Doney~\cite{Alili-Doney}.   These results use special one-dimensional structure of the
 process. For Brownian motion on a half-space, the Martin boundary was
 obtained in the book of Doob~\cite{Doob} by using an explicit form of the
 Green's function. Kurkova
and  Malyshev~\cite{Kurkova-Malyshev} described the  Martin boundary  for random  walks on
$\Z\times\N$ and  on $\Z^2_+$ for which the  only non-zero transitions in  the interior of
the domain are on the nearest neighbors:  $p(z,z \pm e_i) = \mu(\pm e_i)$ with $e_1=(1,0)$
and  $e_2=(0,1)$. For such  random walks,  the jump  generating function  is defined  by $
\varphi(x,y) ~=~  \mu(e_1) x  + \mu(-e_1)x^{-1}  + \mu(e_2)y +  \mu(-e_2)y^{-1} $  and the
equation $  xy(1-\varphi(x,y)) ~=~  0 $ determines  an elliptic  curve ${\bf S}$  which is
homeomorphic to  the torus.   To identify  the Martin boundary,  a functional  equation is
derived for  the generating function  of the Green's  function and the asymptotics  of the
Green's function are  calculated by using the methods of complex  analysis on the elliptic
curve ${\bf S}$.  Such a method seems  to be unlikely to apply in more general situations,
for  higher dimensions  or when  the  jump sizes  are arbitrary,  because the proof is
based on the geometrical
properties of the elliptic curve ${\bf S}$.

\medskip 

In the present paper, we identify the Martin boundary for a random walk $Z_+(t)$ on $\Z^d$
which is killed when leaving  the half-space $\Z^{d-1}\times\N^*$ where
$\N^*=\N\setminus\{0\}$. This is a substochastic 
Markov process on $\Z^{d-1}\times\N^*$ with transition matrix
\[
(p(z,z')=\mu(z'-z); \; z,z'\in \Z^{d-1}\times\N^*) 
\]
where $\mu$  is a  probability measure  on $\Z^d$.  Such  a Markov  process dies  when the
homogeneous random walk $Z(t)$ exits  from $\Z^{d-1}\times\N^*$ and is identical to $Z(t)$
until the first time when $Z(t)\not\in\Z^{d-1}\times\N^*$. The Green function
\[
G_+(z,z') = \sum_n \P_z(Z_+(n)=z') 
\]
is the mean number of visits to the point $z'$ starting from $z$ before hitting
the set $\Z^{d-1}\times(-\N)$. The homogeneous random walk $Z(t)$ is assumed to satisfy the
following conditions 

\smallskip 
\noindent
{\bf (A)} {\em The Markov chain $Z(t)$  is irreducible and has a non zero mean}
\[ 
m ~\dot=~ \sum_{z\in\Z^d} ~z~\mu(z) ~\not=~ 0, 
\]
{\em the last coordinate $Y(t)$ of $Z(t)$ is an aperiodic random walk on $\Z$ 
and the  jump generating function $\varphi$ defined by \eqref{e1-0} 
is finite everywhere on $\R^d$.} 

\smallskip

For such a killed random walk, the limit of the sequence of functions 
\[
G_+(\cdot,z_n)/G_+(z_0,z_n)
\]   
can be identified by using the results of Borovkov~\cite{Borovkov:02} for those sequences
$z_n\in Z^{d-1}\times\N$ for which 
$\lim_{n\to\infty}z_n/|z_n| = (u,v)\in\R^{d-1}\times\R_+$  with $v\not=0$. The most difficult case is when
$\lim_{n\to\infty}z_n/|z_n| \in\R^{d-1}\times\{0\}$, in particular if 
\[
z_n = (x_n,y) \in\Z^{d-1}\times\{y\} 
\]
for all $n\in\N$.  Here, the results of
Borovkov~\cite{Borovkov:02} do not work. The analytical method of Kurkova
and Malyshev do not seem  to apply at
all because one should consider the curve determined
by the equation
\[
\sum_{k=(k_1,\ldots,k_d)\in\Z^d} \mu(k) z_1^{k_1}\cdot \ldots \cdot z_d^{k_d} ~=~ 1
\]
with an infinite number of terms at  the left-hand side. 

While  the  problem of  Martin  boundary  identification of  these  killed  random walks  is
interesting  in its own  right, our methods  also apply to the setting of
Ney and Spitzer~\cite{Ney-Spitzer} and lead to another simple proof of their theorem (see
Section~\ref{sec7} below). The technical results of the present paper are  important  
for the identification of the Martin boundary  of  random walks on  the half-space  with reflected  boundary
conditions  on the hyper-plane $\Z^{d-1}\times\{0\}$, see Ignatiouk~\cite{Ignatiouk:08}.

To formulate our result, it is convenient to introduce the following notations~: 
$ q \to a( q )$ is the inverse mapping of the function~\eqref{e1-1} and, for
$ q \in\R^d\setminus\{0\}$, $a( q ) \dot=~ a( q /| q |)$ is the unique point on the
boundary $\partial D$ of the set $D$ where the normal cone to $D$ contains the vector $ q
$. According to this definition, for $q\not= 0$, 
\[
\sup_{a\in D} ~a\cdot q ~=~ a(q)\cdot q ~>~ a\cdot q, \quad \forall a\in
D\setminus\{a(q)\}
\]
(see Rockafellar~\cite{R}). Define also the half-sphere 
$
{\cal S}^d_+ ~\dot=~ \{ q\in \R^{d-1}\times\R_+ :~ |q| =1\} 
$
with $\R_+ =[0,+\infty[$ and the sets $\partial_+D=\{a\in \partial D: q(a) \in \R^{d-1}\times\R_+\}$ and 
\[
 \partial_0D=\{a\in
\partial D:q(a)\in\R^{d-1}\times\{0\} \} \subset\partial_+ D. 
\]
Denote by $X(t)$ [resp. $Y(t)$ ] the vector of $d-1$ first coordinates [resp. the last
  coordinate ]  of the vector $Z(t)\in\R^d$ and  let $\tau$ be the killing time of the
  process $(Z_+(t))$:
\[
\tau = \inf\{ t~: Y(t) \leq 0\}. 
\]
The variables  $X_+(t)$  and $Y_+(t)$ are defined in a similar way for the vector $Z_+(t)$. 
According to the definition of the process $Z_+(t)$, $(X_+(t),Y_+(t)) = (X(t),Y(t))$ for $t <\tau$ and
$Z(\tau) = (X(\tau),Y(\tau))\in\Z^{d-1}\times(-\N)$ is the point where the process
$Z_+(t)$ dies. For every $a\in\partial_+ D$, the function $h_{a,+}$ on
$\Z^{d-1}\times\N^*$ is defined by    
\[
h_{a,+}(z) ~\dot=~  \begin{cases}
y ~\exp(a\cdot z)  -  \E_{z}\bigl( Y(\tau) ~\exp(a\cdot Z(\tau)), \; \tau < +\infty
\bigr) &\text{if $a \in\partial_0 D$,}\\
\exp(a \cdot z) -  \E_{z}\bigl( \exp(a\cdot Z(\tau)), \; \tau < +\infty\bigr)
&\text{if $a\in\partial_+ D\setminus\partial_0 D$.} 
\end{cases}
\]
The main result of the paper is the following theorem.

\begin{theorem}\label{th1-4} Under the hypotheses (A), the following assertions
  hold : \\
1) The constant
multiples of  the functions $h_{a,+}$ with $a\in\partial_+ D$ are the only minimal harmonic functions of the
Markov process $(Z_+(t))$. \\ 
2) For any $q\in {\cal S}^{d}_+$ and any sequence of points
  $z_n\in\Z^{d-1}\times\N^*$ with $\lim_n|z_n| = \infty$ and $\lim_n z_n/|z_n| = q$, 
\be\label{e1-3}
\lim_{n\to\infty} ~{G_+(z,z_n)}/{G_+(z_0,z_n)} ~=~ h_{a(q),+}(z)/h_{a(q),+}(z_0) 
\ee
for all $z=(x,y)\in\Z^{d-1}\times\N^*$.
\end{theorem}
From this theorem it follows the following statement.
\begin{cor}\label{cor1.1} Under the hypotheses (A), the following assertions hold :\\
1) A sequence of points
$z_n\in\Z^{d-1}\times\N^*$ with $\lim_n |z_n| = 
+\infty$ converge to a point of the Martin boundary for the Markov process $Z_+(t)$   if and
only if $z_n/|z_n| \to q$ for 
some point $q$ on the half-sphere ${\cal S}^d_+$.\\
 2) The full Martin compactification of the half-space 
$\Z^{d-1}\times\N^*$ is homeomorphic to the closure of the set $\left\{ w =
{z}/{(1+|z|)}~:~ z\in\Z^{d-1}\times\N^*\right\}$ in $\R^d$. \\
3) The minimal Martin boundary
coincides with the whole Martin boundary. 
\end{cor}

\medskip 
The proof of Theorem~\ref{th1-4} uses the properties  of Markov-additive processes. A
Markov process $(A(t),M(t))$ on a countable set 
$\Z^d\times E$ is said to be Markov-additive if its transition probabilities are invariant
with respect to the shifts on $x\in\Z^d$ : 
\[
\P_{(x,y)}\bigl((A(t),M(t))=(x',y')\bigr) = \P_{(0,y)}\bigl((A(t),M(t))=(x'-x,y')\bigr),
\]
for all $(x,y),(x',y')\in\Z^d\times E$. For such a Markov process, $A(t)$ is called  an
additive part and $M(t)$ is its Markovian part. With this definition, the Markov process 
$(Z_+(t)=(X_+(t),Y_+(t)))$ is Markov-additive with an additive part $X_+(t)$ on
$\Z^{d-1}$ and a Markovian part $Y_+(t)$ on $E=\N^*$. 

The first assertion of Theorem~\ref{th1-4}  is proved by using the arguments of
Choquet-Deny theory adapted for Markov-additive processes. 
The main steps of the proof of the second assertion are the following. 

Under some general assumptions on the Markov-additive
process $(A(t),M(t))$ on $\Z^d\times E$, it is proved that 
the  Green's function  
\[
{\cal G}((k,y),(l,y')) ~=~ \sum_n \P_{(k,y)}\bigl((A(n),M(n))=(l,y')\bigr)
\]
satisfies  the following  property  :  if for  a  sequence $z_n=(x_n,y_n)\in\Z^d\times  E$
converging to infinity, the inequality 
\[
\lim_{n\to\infty} \frac{1}{|z_n|} \log {\cal G}((0,y),z_n) \geq 0 
\]
holds for some $y\in E$ then also 
\[
\lim_{n\to\infty} {\cal G}((x+ x',y),z_n)/{\cal G}((x + x'',y),z_n) = 1, \quad \quad \forall 
x,x'\in\Z^d, \; y\in E
\]
for all $(x,y)\in\Z^d\times E$ and all those $x',x''\in\Z^ð$ for which 
\[
\inf_{y\in E} ~\min\Bigl\{\P_{(0,y)}(A(k) = x', \; M(k) = y), \; \P_{(0,y)}(A(k) = x'', \;
M(k) = y) \Bigr\} ~>~ 0. 
\]
This is an extension of an intermediate  technical result of
Foley and McDonald~\cite{Foley-McDonald}) 
obtained for a different purpose in a more restricted context for Markov-additive processes with a
one-dimensional additive part and for a 
sequence of the form $z_n=(n,y)$.  

For the random walk $(Z_+(t))$, this result implies that 
\be\label{e1-4}
\lim_{n\to\infty} G_+((x + \hat{k}x',y),z_n)/G_+((x,y),z_n) ~=~ 1, \quad \quad \forall x,x'\in\Z^{d-1}, \; y\in \N
\ee
whenever 
\be\label{e1-5}
\lim_{n\to\infty} \frac{1}{|z_n|} \log G_+(z,z_n) ~=~ 0 \quad \quad \text{for some} \quad 
  z\in\Z^{d-1}\times\N^*
\ee
where $\hat{k}$ denotes the period of the random walk
$(Z(t))$. 

The second important tool of our proof is the use of large deviation techniques.  With the
aid of 
Mogulskii's theorem~(see Dembo and Zeitouni~\cite{D-Z}), we show that the family
of scaled processes $Z^\eps_+(t) = \eps Z_+([t/\eps])$ satisfies sample path large deviation
principle. 
The inequality \eqref{e1-5} is obtained for any 
$z\in\Z^{d-1}\times\N^*$ and any  sequence of points
$z_n\in\Z^{d-1}\times\N^*$ when 
\[
\lim_{n\to\infty} z_n/|z_n| ~=~  q(0)\in\partial_+ D 
\]
by using the lower large deviation bound. 

Finally, we use the results of Doney~\cite{Doney:02} and the integral 
representation of the harmonic functions to show
that for every $a=(\a,\beta)\in\partial 
_+ D$ with $\a\in\R^{d-1}$ and $\beta\in\R$, the constant multiples of the function
$h_{a,+}$ are the only non-negative harmonic 
functions satisfying the equality 
\[
h(x +\hat{k} x',y)~=~ \exp(\hat {k} \a\cdot x') h(x,y), \quad \quad \forall x,x'\in\Z^{d-1}, \; y\in
\N. 
\]
The last result and the equality ~\eqref{e1-4} are used to get the convergence \eqref{e1-3} when $q=q(0)\in
\R^{d-1}\times\R_+$.  An exponential change of measure with a parameter $a(q)$ 
extends this result for an arbitrary $q\in \cal{S}_+^d$.

Before proving our results, an example where the functions $h_{a,+}$ have an explicit form is given.

\section{An explicit representation of the functions $h_{a,+}$ } 
In this section, we calculate explicitly the functions $h_{a,+}$ for $a\in\partial_+ D$ in
a particular   case  when  $\mu(z)  =   0$  for  all 
$z=(x,y)\in\Z^{d-1}\times\Z$ with $y<-1$.  In this  classical situation, the last vertical 
component $(Y(t))$ of  the random walk $(Z(t))$ is  left continuous, i.e. the  value of a  downward jump is
$-1$.

Let  $\pi : \R^d \to \R^d$  be the orthogonal
projection to  the hyper-plane  $\R^{d-1}\times\{0\}$. For  a point  $a\in\partial_+ D$
with $q(a)\in S^d_+$ we  denote by $\ol{a}$ the unique point on  the boundary $\partial D$
with   $q(\ol{a})\in  S^d\cap(\R^d\times]-\infty,0])$   for   which  $\pi(a)=\pi(\ol{a})$.
Clearly, $a=\ol{a}$ if and only if $a\in\partial_0 D$.

\begin{prop}\label{pr2-1} Suppose that $\mu(z) = 0$  for all $ z=(x,y)\in\Z^d$ with  $y< -
  1$, then under the hypotheses (A), for every $a\in \partial_+ D$,
\be\label{e2-1}
h_{a,+}(z) ~=~ \begin{cases} \exp(a\cdot z) - \exp(\ol{a}\cdot z) &\text{if
    $a\not\in\partial_0 D$,}\\
y ~\exp(a\cdot z) &\text{if
    $a\in\partial_0 D$,} \quad \quad \; \forall z=(x,y)\in\Z^d\times\N^*.
\end{cases}
\ee
\end{prop}
\begin{proof} Indeed, if $\mu(z) ~=~ 0$  for all  $z = (x,y)\in\Z^{d-1}\times\Z$  with  $y
    < - 1$ then almost surely $Y(\tau) = 0$ and consequently, for every
    $a=(\alpha,\beta)\in\partial_+ D$ with $\alpha\in\R^{d-1}$ and $\beta\in \R$, 
\[
h_{a,+}(z) ~=~ \begin{cases}  \exp(a\cdot z) - \E_{z}\left(
  \exp\bigl(\alpha\cdot X(\tau)\bigr); \; \tau < +\infty\right) &\text{if
   \; $ a\in\partial_+ D\setminus \partial_0 D$,}\\
y \, \exp(a\cdot z) &\text{if
   \; $ a\in\partial_0 D$.}
\end{cases}
\]
Furthermore, for $a\in\partial_+ D\setminus\partial_0 D$, the quantity  
\[
\exp(-\ol{a} \cdot z) \E_{z}\left(\exp\bigl(\alpha\cdot X(\tau)\bigr); \; \tau < +\infty\right) ~=~
\E_{z}\left(\exp\bigl(\ol{a}\cdot (Z(\tau) -z)\bigr), \; \tau < +\infty\right) 
\]
is equal to the probability that the twisted homogeneous random walk $\ol{Z}(t)$ on $\Z^{d}$
with transition probabilities 
\[
p_{\ol{a}}\bigl(z,z'\bigr) ~=~ \exp\bigl(\ol{a}\cdot(z'-z)\bigr) \mu(z'-z),
\]  
starting from $z$ ever hits the boundary hyper-plane $\Z^{d-1}\times\{0\}$. Since the  last
coordinate of the mean  
\[
\E_0(\ol{Z}(1)) ~=~
  \sum_{z\in\Z^{d-1}\times\Z} z \,\mu(z) 
  \exp(\ol{a}\cdot z) ~=~ \nabla \varphi(\ol{a}) 
\]
is negative then for every starting point 
$z\in\Z^{d-1}\times\N^*$, this hitting probability is equal to $1$ and consequently, 
\eqref{e2-1} holds. 
\end{proof}

\section{General properties of Markov-additive processes}\label{sec3}
In this section we describe minimal harmonic function of Markov-additive processes and
prove the ratio limit theorem. Before to prove these results we recall the definition of
{\it Markov-additive processes}, the corresponding {\it Feynman-Kac transform}  and its
{\it spectral radius}. 

A Markov process 
$(A(t),M(t))$ on a countable set $\Z^{d}\times E$ with transition probabilities  
  $p\bigl((x,y),(x',y')\bigr)$ is called {\em Markov-additive} if 
\[
p\bigl((x,y),(x',y')\bigr) ~=~ p\bigl((0,y),(x'-x,y')\bigr)
\]
for all $x,x'\in\Z^{d}$, $y,y'\in E$. The first component $A(t)$ is an {\em additive} part of the process
$(A(t),M(t))$,  and $M(t)$ is its {\em Markovian part}. The Markovian part $M(t)$ is a Markov
chain on $E$ with transition  probabilities 
\[
P(y,y') = \sum_{x\in\Z^{d}} p\bigl((0,y),(x,y')\bigr). 
\]
For $\a\in\R^d$, an infinite matrix ${\cal P}(\alpha) = \bigl( {\cal P}(\alpha, y,y'), \; y,y'\in E\bigr)$
with 
\[
{\cal P}(\alpha, y,y') ~=~ \E_{(0,y)}\bigl( \exp(\alpha\cdot A(1)) ; \, M(1) = y'\bigr)
\]
is called {\it Feynman-Kac transform} of the Markov-additive process $(A(t),M(t))$. The
$n$-th iterate  ${\cal P}^{n}(\alpha) = ( {\cal P}^{(n)}(\alpha, y,y'), \; y,y'\in
E\bigr)$ of the matrix  ${\cal P}(\alpha)$ is given by 
\[
{\cal P}^{(n)}(\alpha, y,y') ~=~ \E_{(0,y)}\bigl(
\exp(\alpha\cdot A(n)) ; \, M(n) = y'\bigr), \quad y,y'\in E. 
\] 
If the Markovian part $(M(t))$ is irreducible then the matrix ${\cal P}(\alpha)$ is
also irreducible. In this case, for every $\lambda\in\R$, the series 
\be\label{e3-1}
G_{\lambda}(\alpha, y,y') ~\dot=~ \sum_{n=0}^\infty e^{- \lambda n} {\cal P}^{(n)}(\alpha,
y,y')  
\ee
converge or diverge simultaneously for all $y,y'\in E$, and   
the limit
\[
\lambda(\alpha) ~=~ \limsup_n \frac{1}{n} \log {\cal P}^{(n)}(\alpha, y,y') 
\]
does not depend on $y,y'\in E$ (see~\cite{Seneta}). The quantity $e^{\lambda(\alpha)}$ is usually called
  {\it spectral radius} and $e^{-\lambda(\alpha)}$ is the {\it convergence parameter} of the
  transform matrix ${\cal P}(\alpha)$, it is a common radius of
  convergence of the series \eqref{e3-1}. According to this definition, $\lambda(\alpha) = +
  \infty$ if ${\cal P}(\alpha,y,y') = +\infty$ for some $y,y'\in E$. 

\subsection{Minimal harmonic functions}

Recall that a non-negative function $h$ on $\Z^{d}\times E$ is {\em harmonic} for the
Markov process $(A(t),M(t))$ if $\E_z(h(A(t),M(t))) = h(z)$ for all $z\in\Z^{d}\times E$. A 
non-zero harmonic function $h\geq 0$ is
  {\em minimal} if for any non-zero harmonic function $h'\geq 0$, the inequality $h'\leq h$
  implies that $h'=c h$ with some constant $c > 0$. 

\begin{prop}\label{pr3-0} Suppose that the Markovian part $(M(t))$ is irreducible and that 
  for every $x\in\Z^{d}$, there are $n\in\N$ and $\theta>0$ such that
  $p^{(n)}\bigl((0,y),(x,y)\bigr) \geq \theta$ for all $y\in E$. Then
  every non-zero minimal harmonic 
  function $h$ of the Markov process $(A(t),M(t))$ is of the form 
  \be\label{e3-2}
h(x,y) ~=~ e^{\alpha \cdot x} h(0,y) ~>~ 0, \quad \quad \forall (x,y)\in \Z^{d}\times E
\ee 
with some  $\alpha\in\R^{d}$  satisfying
  the inequality $\lambda(\alpha) \leq 0$.  
\end{prop}
\begin{proof} The proof of this proposition uses the arguments similar to that of Choquet and Deny
  theorem  (see Woess~\cite{Woess}). 

Let $h(x,y) \geq 0$ be a 
  harmonic function. Then for a unit vector $e_i\in\Z^{d}$, the function $h_i(x,y) = h(x + e_i,
  y)$ is also harmonic and under the hypotheses of our proposition, there exist $n_i\in\N$
  and $\theta >0$ such that  
\begin{align*}
h(x,y) &= \sum_{(x',y')} p^{(n_i)}\bigl((x,y),(x',y')\bigr) h(x',y') \geq p^{(n_i)}\bigl((x,y),
(x+e_i,y)\bigr) h(x+ e_i,y) \\&\geq~ \theta h_i(x,y)
\end{align*}
for all $(x,y)\in\Z^{d}\times E$. If the harmonic function $h$ is minimal, the last
inequality implies that $h_i = c_i h$ for some $c_i > 0$. Using the equality $h= c_i h_i$ for
every unit vector $e_i\in\Z^{d}$, 
$i=1,\ldots, d$, and letting $\alpha_i = \ln c_i$ we
obtain~\eqref{e3-2}. Moreover, since the function $h$ is harmonic and since  
~\eqref{e3-2} holds, 
\[
h(0,y) = \sum_{y'\in E} {\cal P}^{(n)}(\alpha, y,y') h(0,y') \geq {\cal P}^{(n)}(\alpha, y,y') h(0,y') 
\]
for all $y,y'\in E$ and for all $n\in\N$. Under the hypotheses of our proposition,  the
matrix ${\cal P}(\alpha)$ is  irreducible and hence, the last inequality implies that
$h(0,y) > 0$ for all $y\in E$ and 
\[
\lambda(\alpha) = \limsup_n \frac{1}{n} \log {\cal P}^{(n)}(\alpha, y,y')  \leq 0  
\]
whenever  $h\not\equiv 0$. 
Proposition~\ref{pr3-0} is therefore proved. 
\end{proof}

\subsection{Ratio limit theorem}

Throughout this section, ${\cal Z}(t)=(A(t),M(t))$ denotes a Markov-additive process on
$\Z^d\times E$ with an additive part $(A(t))$ on $\Z^d$ and a Markovian part $(M(t))$ on $E\subset\Z^k$. 
Green's function of the Markov process ${\cal Z}(t)$ is denoted by ${\cal
  G}(z,z')$. The 
 following conditions are assumed to be satisfied. 

\medskip 
\noindent
{\bf (A1)} {\em The Markov-additive process ${\cal Z}(t)=(A(t),M(t))$ is irreducible. }

\noindent 
{\bf (A2)} {\em The function 
\be\label{e3-3}
\hat\varphi(a) ~=~ \sup_{z\in \Z^d\times E} 
~\E_{z}\Bigl( \exp\bigl(a\cdot ({\cal Z}(1)-z)\bigr)\Bigr)   
\ee
is finite in a neighborhood of zero. }

\medskip

Remark that the Markov-additive process ${\cal Z}(t)=(A(t),M(t))$ is not necessarily
stochastic~: in some points $z = (x,y)\in\Z^d\times E$, 
the transition matrix can be strictly substochastic.  
\medskip

\begin{prop}\label{pr3-1} Suppose that  the Markov-additive process ${\cal
    Z}(t)=(A(t),M(t))$ is transient and satisfies 
  the hypotheses (A1), (A2) 
and let a sequence of points $z_n\in\Z^d\times E$
  be such that  $|z_n|\to\infty$ as $n\to\infty$ and for some $z_0\in\Z^d\times E$, 
\be\label{e3-4}
\liminf_{n\to\infty}~\frac{1}{|z_n|} \log {\cal G}\bigl(z_0, z_n\bigr)  ~\geq~ 0. 
\ee
Then  
\[
\lim_{n\to\infty} {\cal G}(z+w,z_n)/{\cal G}(z+w',z_n) ~=~ 1. 
\]
for all $z\in\Z^d\times E$ and for all those $w,w'\in\Z^d\times\{0\}$ for which there is
$n > 0$ such that 
\[ 
 \inf_{z\in \Z^d\times E} ~\min\left\{p^{(n)}\bigl(z, z + w\bigr), p^{(n)}\bigl(z,
z + w' \bigr)\right\} ~>~ 0. 
\] 
\end{prop}

This statement was initially obtained by Foley and McDonald~\cite{Foley-McDonald} for
Markov additive processes with an additive part on $\Z$ and for a sequence of the form
$z_n=(n,y)$ with a given $y\in E$. In general case, the proof of this proposition uses
essentially the same ideas as in ~\cite{Foley-McDonald}, we give this proof in 
 section~\ref{secA-1}.

Consider now  the following more restrictive conditions.  

\medskip 
\noindent
{\bf (A1')} {\em (Communication condition) There exist 
$\theta >0$ and $C>0$ such that for any $z\not=z'$, $z,z'\in\Z^d\times E$ there is a sequence of
  points $z_0, z_1, \ldots,z_n\in\Z^d\times E$  with $z_0=z$, $z_n=z'$ and 
    $n\leq C|z'-z|$ such that  
\[
 |z_i-z_{i-1}| \leq C \quad \text{ and  } \quad \P_{z_{i-1}}({\cal Z}(1) = z_i) \geq \theta, \quad \quad \forall \;
 i=1,\ldots,n. 
\]}
{\bf \!(A2')} {\em For every $z\in\Z^d\times E$, the function 
\[
\varphi_z(a) ~=~ 
~\E_{z}\Bigl( \exp\bigl(a\cdot ({\cal Z}(1)-z)\bigr)\Bigr)   
\]
is finite everywhere on $\R^{d+k}$ and  the function $\hat\varphi$ defined by \eqref{e3-3} 
is finite in a neighborhood of zero.}

\medskip

If the assumption (A1') is satisfied then there is a bounded function
$n_0 : E \to \N^*$ such that for any $z=(x,y)\in\Z^d\times E$, 
\[
p^{(n_0(y))}\bigl((x,y), (x,y)\bigr) \geq \theta^{n_0(y)} > 0
\]
and hence, there is $k  \in\N^*$  (for instance, $k = n!$ with $n = \max_y n_0(y)$) such 
that 
\[
 p^{(k)}\bigl(z, z\bigr) \geq  \theta^{k}, \quad \quad \quad \forall z\in\Z^d\times E.
\]
We denote by ${\cal K}$ the set of all integers $k >0$ for
which 
\be\label{e3-5}
\inf_{z\in Z^d\times E} p^{(k)}(z,z) ~>~ 0 
\ee 
The greatest common divisor of the set ${\cal K}$ is denoted by $\hat{k}$ and the
following condition is assumed to be satisfied. 

\medskip
\noindent
{\bf (A3)} {\em Up to multiplication by constants, there is a unique 
 positive harmonic function  $h : \Z^d\times E\to\R_+$ satisfying the equality  
$h(z + \hat{k}w) = h(z)$ for all $z\in\Z^{d}\times E$ and $w\in
 \Z^d\times\{0\}$.}
\medskip 

When
the  Markov-additive 
process ${\cal Z}(t)=(A(t),M(t))$ is stochastic,  the last assumption means that  the only
positive harmonic functions  $h : \Z^d\times 
E\to\R_+$ satisfying the equality   
$h(z+\hat{k}w) = h(z)$ for all $z\in\Z^{d}\times E$ and $w\in
 \Z^d\times\{0\}$ are constant. 

Under these additional assumptions, using Proposition~\ref{pr3-1} we obtain the following statement. 

\begin{prop}\label{pr3-2} Suppose that  the  
Markov-additive process ${\cal Z}(t)=(A(t),M(t))$ is transient and satisfies the hypotheses
 (A1') and (A2), and let a sequence of points $z_n\in\Z^d\times E$ with $\lim_n |z_n| =
 \infty$ satisfy the inequality \eqref{e3-4} for some $z_0\in\Z^d\times E$. Then 
\be\label{e3-6}
\lim_{n\to\infty} {\cal G}(z + \hat{k} w,z_n)/{\cal G}(z,z_n) ~=~ 1,  \quad \quad \quad  \forall z\in\Z^d\times
E, \; w\in\Z^d\times\{0\}. 
\ee
Moreover, if  the
conditions (A2') and (A3) are also satisfied, then 
\be\label{e3-7}
\lim_{n\to\infty} {\cal G}(z,z_n)/{\cal G}(z',z_n) ~=~ h(z)/h(z'), \quad \quad \forall z,z'\in\Z^d\times E.
\ee
\end{prop}
\begin{proof} Let us show that 
\be\label{e3-8}
\lim_{n\to\infty} {\cal G}(z + k w,z_n)/{\cal G}(z,z_n) ~=~ 1
\ee
for all $k\in{\cal K}, z\in\Z^d\times E$ and $w\in\Z^d\times\{0\}$. Indeed, according to
the definition of the set ${\cal K}$, if $k\in{\cal
  K}$ then there is $\eps > 0$ such that 
\be\label{e3-9}
p^{(k)}(z,z) ~\geq~ \eps, \quad \quad \quad \forall z\in\Z^d\times E. 
\ee
Furthermore, because of the assumption (A1'), for every  $w\in 
\Z^d\times\{0\}$, there is a bounded positive function $n_w(\cdot) : E 
\to \N^*$ with $n_w \dot= \sup_{y\in E} n_w(y) \leq C|w|$ such that for any
$(x,y)\in\Z^d\times E$, 
\[
 p^{(n_w(y))}\bigl((x,y), (x,y)+ w\bigr) ~=~ p^{(n_w(y))}\bigl((0,y), (0,y)+ w\bigr)
 ~\geq~ \theta^{n_w(y)} ~\geq~ \theta^{n_w}. 
\]
From this it follows that  
\be\label{e3-10}
p^{n_w(y)}\bigl((x,y)+ (j-1)w, (x,y) + jw\bigr) ~\geq~ \theta^{n_w} 
\ee 
for all $(x,y)\in\Z^d\times E$, $w\in\Z^d\times\{0\}$ and $j\in\N^*$ (to get the last 
inequality we replace in the above inequality $(x,y)$ by $(x,y) +
(j-1)w$).  Using the  inequalities \eqref{e3-9} and \eqref{e3-10} we get  
\begin{align*}
p^{(k n_w)}(z,z + k w) &~\geq~ \left(p^{(k)}(z,z)\right)^{n_w-n_w(y)}
\prod_{j=1}^{k} p^{n_w(y)}\bigl(z+ (j-1)w, z +
jw\bigr) \\ &~\geq~ \eps^{n_w-n_w(y)} \theta^{k n_w} ~\geq~
\eps^{n_w} \theta^{k n_w}   
\end{align*}
for all $(x,y)\in\Z^d\times E$, $w\in\Z^d\times\{0\}$ and consequently,   
\[
\inf_{z\in \Z^d\times E} ~\min\left\{p^{(k n_w)}\bigl(z, z\bigr), p^{(k n_w)}\bigl(z,
z + \hat{k}w\bigr)\right\} ~\geq~ \min\left\{ \eps^{n_w}, \eps^{n_w}\theta^{k n_w}\right\}
~>~ 0 
\]
for every $w\in\Z^d\times\{0\}$.  By  Proposition~\ref{pr3-1}, under the hypotheses (A1')
and (A2), the last inequality proves \eqref{e3-8}. 

\medskip

Consider now the subgroup $\langle {\cal K} \rangle$ of $\Z$ generated by the set ${\cal
  K}$ and let us notice that   for every  $k \in{\cal K}$, we have also 
\[
\lim_{n\to\infty} {\cal G}(z - kw ,z_n)/{\cal G}(z, z_n) ~=~ 1, \quad \quad \quad  \forall z\in\Z^d\times
E, \; w\in\Z^d\times\{0\}  
\]
To get this relation it is sufficient to replace  $z$ by $z -
k w$ in  \eqref{e3-8}. Hence, \eqref{e3-8} holds for every $k \in\langle {\cal K} \rangle$ and
in particular for $k = \hat{k}$ because $\hat{k}\in\langle{\cal K}\rangle$ (see Lemma~A.1
of Seneta~\cite{Seneta}). The first assertion of Proposition~\ref{pr3-2} is therefore proved. 

\medskip 

Suppose now that the conditions (A2') and (A3) are also satisfied. 
Because of the assumption (A1'), for any $z,z'\in\Z^d\times E$, the
  probability that  the Markov process 
${\cal Z}(t)$ starting 
  at $z$ ever hits $z'$ is greater than $\theta^{C |z'-z|}$ which implies that 
\be\label{e3-11}
\theta^{C |z-z_0|}   ~\leq~ {\cal G}(z,z_n)/{\cal G}(z_0,z_n) ~\leq~ \theta^{-C |z-z_0|}  
\ee 
for all $z\in\Z^d\times E$.  Suppose now that for a  
subsequence $(z_{n_k})$, the sequence of functions 
\[
 {\cal G}(z,z_{n_k})/ {\cal G}(z_0,z_{n_k}), \quad z\in\Z^d\times E
\]
converge pointwise in $\Z^d\times E$. Then from the first inequality of \eqref{e3-11} it
follows that the function  
\[
\tilde{h}(z) ~\dot=~ \lim_{k\to\infty}  {\cal G}(z,z_{n_k})/ {\cal G}(z_0,z_{n_k})
\]
is strictly positive. By dominated convergence theorem, the second inequality of
\eqref{e3-11} implies that the function $\tilde{h}$ is harmonic for $({\cal Z}(t))$ (recall
that under the hypotheses (A2'), the exponential functions
are integrable with respect to the probability measure $p(z,\cdot)$ for every
$z\in\Z^d\times E$).  Moreover, the equality \eqref{e3-6} shows that  
   $\tilde{h}(z) = \tilde{h}(z + \hat{k}w)$ 
for all $z\in\Z^d\times E$ and $w\in\Z^d\times\{0\}$.  Since the constant multiples of the
function $h$ are the only harmonic functions satisfying this equality, we conclude that
   $\tilde{h} = c h$ with some $c>0$.  Moreover, since $\tilde{h}(z_0) = 1$ then 
$\tilde{h}(z) ~=~ h(z)/h(z_0)$ and consequently, 
\[
\lim_{k\to\infty}  {\cal G}(z,z_{n_k})/ {\cal G}(z_0,z_{n_k}) ~=~ \tilde{h}(z) ~=~ h(z)/h(z_0), \quad \quad \forall
\;   z\in\Z^d\times E
\]
for every subsequence $n_k$ for which these limits exist. From this it follows that 
\[
\lim_{n\to\infty}  {\cal G}(z,z_{n})/ {\cal G}(z_0,z_{n}) ~=~ \tilde{h}(z) ~=~ h(z)/h(z_0), \quad \quad \forall
\;   z\in\Z^d\times E
\]
because for every $z\in\Z^d\times E$, the sequence $ {\cal G}(z,z_n)/ {\cal G}(z_0,z_n)$ is bounded. 
\end{proof}

\section{Sample path large deviations}\label{sec4}    

\subsection{General statement of SPLD principle}
Let $D([0,T],\R^{d})$ denote the set of all right continuous functions with left
limits from $[0,T]$ to $\R^{d}$  endowed with Skorohod metric
(see Billingsley~\cite{Billingsley}).

\smallskip
\noindent
\begin{defi} 
1) A mapping $I_{[0,T]}:~D([0,T],\R^{d})\to
[0,+\infty]$ is a good rate function on $D([0,T],\R^{d})$ if for any $c\geq 0$ and 
 any compact set $V\subset \R^{d}$, the set
\[
\{ \phi \in D([0,T],\R^{d}): ~\phi(0)\in V \; \mbox{
and } \; I_{[0,T]}(\phi) \leq c \}
\]
is compact in $D([0,T],\R^{d})$.  According to this definition, a good
rate function is lower semi-continuous. 

2) Let $(Z(t))$ be a Markov process on  $E\subset \Z^d$ and let 
$Z^\eps(t) =\eps Z([t/\eps])$ for $\eps>0$.
When $\eps\to 0$, the family of scaled processes $(Z^\eps(t) =\eps Z([t/\eps]),
\,t\in[0,T])$,  is said to
satisfy {\it sample path large deviation principle} in $D([0,T], \R^{d})$ with a rate function
$I_{[0,T]}$ if for any $z\in\R^{d}$ 
\begin{equation}\label{e4-1}
\lim_{\delta\to 0} \;\liminf_{\eps\to 0} \; \inf_{z'\in \eps E : |z'-z|<\delta} \eps
\log\P_{[z'/\eps]}\left( Z^\eps(\cdot)\in {\cal 
O}\right) \geq -\inf_{\phi\in{\cal O}:\phi(0)=z} I_{[0,T]}(\phi), 
\end{equation}
for every open set ${\cal
O}\subset D([0,T],\R^{d})$,
and
\begin{equation}\label{e4-2}
\lim_{\delta\to 0} \;\limsup_{\eps\to 0} \; \sup_{z' \in \eps E : |z'-z|<\delta}
\eps\log\P_{[z'/\eps]}\left( Z^\eps(\cdot)\in 
F\right) \leq -\inf_{\phi\in F:\phi(0)=z} I_{[0,T]}(\phi).
\end{equation}
 for every closed set $F\subset   D([0,T],\R^{d})$. 
\end{defi} 
$\P_{[z/\eps]}$ denotes here and throughout the distribution of the Markov process $(Z(t))$
corresponding to the initial state $Z(0)=[z/\eps]$ where $[z/\eps]$ is the nearest lattice
point to $z/\eps$ in $E\subset\Z^d$.  
We refer to sample path large deviation
principle as SPLD principle. Inequalities (\ref{e4-1}) and (\ref{e4-2}) are
referred as lower and upper SPLD bounds respectively.

\subsection{SPLD properties of scaled processes $Z_+^\eps(t) = \eps
    Z_+([t/\eps]) $} Under the 
    hypotheses (A), the jump generating function $\varphi$ of \eqref{e1-0}  
is finite everywhere on $\R^d$ and hence, by Mogulskii's theorem~(see~\cite{D-Z}), 
the family of scaled processes $Z^\eps(t)=\eps Z([t/\eps])$ satisfies SPLD principle in
$D([0,T],\R^{d})$ with a 
good rate function 
\[
I_{[0,T]}(\phi) ~=~ \begin{cases} \int_0^T (\log \varphi)^*(\dot\phi(t)) \, dt, &\text{ if
    $\phi$ is absolutely continuous,}\\
+\infty &\text{ otherwise.}
\end{cases}
\]
Recall that a continuous function $\phi : [0,T]\to\R^d$ is  absolutely continuous if and
only if it is differentiable almost everywhere 
in $[0,T]$ with $\dot\phi \in L^1([0,T])$ and for every $t\in[0,T]$,
\[
\phi(t) = \phi(0) + \int_0^t \dot\phi(s) \, ds 
\]
(see Rudin~\cite{Rudin}). 
The convex conjugate  $(\log \varphi)^*$ of the function $\log \varphi$ is defined by  
\[
(\log \varphi)^*(v) ~\dot=~ \sup_{a\in\R^{d}} \Bigl(a\cdot v - \log
\varphi(a)\Bigr). 
\]
Under the hypotheses (A), $
(\log \varphi )^*(v) ~=~ a\cdot v - \log
\varphi(a)$ whenever $v ~=~ \nabla \varphi(a)$
because the function $(\log \varphi)$ is convex and differentiable
everywhere in $\R^d$ (see Lemma~2.2.31 of the book of Dembo and Zeitouni~\cite{D-Z}). In
particular, according to the definition of the function $q\to a(q)$, 
\be\label{e4-4}
(\log \varphi )^*\Bigl(\nabla \varphi(a(q))  \Bigr) ~=~ a(q) \cdot \nabla \varphi(a(q)) ~=~
|\nabla \varphi(a(q))| a(q) \cdot q  
\ee
because $\varphi(a(q)) = 1$. For a linear function $\phi(t) = vt$ with $v = \nabla
\varphi(a(q))$ and for $T = 1/|\nabla \varphi(a(q))|$ we have therefore $\phi(T) = q$ and $
I_{[0,T]}(\phi) ~=~ a\cdot q$. 

\medskip

To prove SPLD principle for the family of scaled processes $Z_+^\eps(t) = \eps
Z_+([t/\eps])$ we use the following lemma.

\begin{lemma}\label{lem1-4} If the random walk $(Z(t))$ is irreducible then the
Markov process $(Z_+(t))$ satisfies the communication condition (A1')~: there are 
$\theta >0$ and $C>0$ such that for any $z\not=z', \; z,z'\in\Z^{d-1}\times\N^*$ 
there exists a sequence of points $z_0, z_1, \ldots,z_n\in\Z^{d-1}\times\N^*$  with $z_0=z$, $z_n=z'$ and 
    $n\leq C|z'-z|$ such that  
\[
|z_i-z_{i-1}| \leq C \quad \text{ and } \quad 
 \P_{z_{i-1}}(Z_+(1) = z_i) = \mu(z_i-z_{i-1}) \geq \theta,  \quad \forall \;
 i=1,\ldots,n. 
\] 
\end{lemma}
\begin{proof} Indeed, suppose that the random walk $(Z(t))$ is irreducible. Then for every unit vector
 $e\in\Z^d$, there is a sequence of points $u_1^e, \ldots,u_{n_e}^e\in\Z^{d}$ 
 such that $\mu(u_{i}^e)  > 0 $ for every $i = 1,\ldots,n_e$ 
and 
\[
u_1^e + \ldots + u_{n_e}^e~=~ e. 
\]
For every $z\in\Z^d$ there are therefore 
$
u_1,\ldots,u_n \in \cup_{e\in\Z^d : |e| = 1} \{
u_1^e, \ldots, u_{n_e}^e\}
$
with 
\[
n ~\leq~ \sum_{e\in\Z^d : |e| = 1}  n_e ~z\cdot e 
\]
such that $
z = u_1 + \ldots + u_n$. 
This proves that the Markov process $(Z(t))$ satisfies  the communication condition
(A1') with 
\[
\theta ~=~ \min_{e\in\Z^d : |e| = 1} ~\min_{1\leq i \leq n_e} \mu(u_{i}^e) 
\]
and 
\[
C = \max\left\{ \sum_{e\in\Z^d : |e| = 1}  n_e,  \max_{e\in\Z^d : |e| = 1} ~\max_{1\leq i
  \leq n_e}  |u_i^e| \right\}.
\]
Hence, for any $z,z'\in\Z^{d-1}\times\N^*$, if $z\not=z'$ then there is a sequence of
  points $z_0, z_1, \ldots,z_n\in\Z^d$  with  $z_0=z$ and $z_n=z'$ such that 
    $n\leq C|z'-z|$,  
\[
|z_i-z_{i-1}| ~\leq~ C \quad \text{ and } \quad 
\mu(z_i-z_{i-1}) ~\geq~ \theta ~>~ 0  \quad \forall \;
 i=1,\ldots,n. 
\]   
Moreover, without any restriction of generality we can suppose that there is $1\leq k\leq
n$ such that $u_i ~\dot=~ z_i-z_{i-1}\in\Z^{d-1}\times\N$ for all $i < k$ and $u_i ~\dot=~ 
z_i-z_{i-1}\in\Z^{d-1}\times(-\N)$ for $i > 
k$. Then, $z_i ~\dot=~ z + u_1 + \ldots + u_i \in \Z^{d-1}\times\N^*$ for all
$i=1,\ldots,n$ and consequently 
\[
\P_{z_{i-1}}(Z_+(1) = z_i) ~=~  \mu(z_i-z_{i-1})  ~\geq~ \theta, \quad \quad \forall \;
 i=1,\ldots,n. 
\]
\end{proof}

\begin{prop}\label{pr4-1} If the random walk $(Z(t))$ is irreducible and the jump generating function
  $\varphi$ of \eqref{e1-0} is finite everywhere on $\R^d$ then the sequence of scaled processes
  $Z^\eps_+(t)=\eps Z_+([t/\eps])$ satisfies SPLD principle in 
$D([0,T],\R^{d})$ with the 
good rate function 
\[
I_{[0,T]}^+(\phi) ~=~ \begin{cases} \int_0^T (\log \varphi)^*(\dot\phi(t)) \, dt, &\text{ if
    $\phi$ is absolutely continuous and}\\
&\text{ $\phi(t)\in\R^{d-1}\times\R_+$ for all $t\in[0,T]$,}\\
+\infty &\text{ otherwise.}
\end{cases}
\]
\end{prop}
\begin{proof}
Indeed, for any $c\geq 0$ and for any compact set
$V\subset \R^{d}$, the 
    set 
\begin{multline*}
\{\phi  : I_{[0,T]}^+(\phi) \leq c, \, \phi(0)\in V\} = \{\phi : I_{[0,T]}(\phi) \leq c,
\, \phi(0)\in V\} ~\cap \\ \{\phi : \phi(t)\in\R^{d-1}\times\R_+, \, \forall t\in[0,T]\} 
\end{multline*}
is compact in $D([0,T],\R^{d})$ because $I_{[0,T]}$ is a good rate function on
$D([0,T],\R^{d})$ and the set $ \{\phi : \phi(t)\in\R^{d-1}\times\R_+,\, \forall t \in[0,T] \}$ is closed in
$D([0,T],\R^{d})$. This proves that the mapping $I_{[0,T]}^+ : D([0,T],\R^{d}) \to
[0,+\infty]$ is a good rate function.   

Furthermore, for any closed set $F\subset
D([0,T],\R^{d})$, the set $F_+ ~\dot=~ \{\phi \in F : ~\phi(t)\in\R^{d-1}\times\R_+,\,
\forall t \in[0,T]\}$ is also closed and 
consequently, for any $z\in\R^{d-1}\times\R_+$, using Mogulskii's theorem we obtain 
\begin{align*}
\lim_{\delta\to 0} \;\limsup_{\eps\to 0} \; \sup_{z' 
  : |z'-z|<\delta} \eps &\log\P_{[z'/\eps]}\left( Z^\eps_+(\cdot)\in
F\right) \\&= \lim_{\delta\to 0} \;\limsup_{\eps\to 0} \; \sup_{z'
  : |z'-z|<\delta} \eps \log\P_{[z'/\eps]}\left( Z^\eps_+(\cdot)\in
F_+\right) \\ &\leq \lim_{\delta\to 0} \;\limsup_{\eps\to 0} \; \sup_{z'
  : |z'-z|<\delta} \eps \log\P_{[z'/\eps]}\left( Z^\eps(\cdot)\in
F_+\right)\\
&\leq -\inf_{\phi\in F_+:\phi(0)=z} I_{[0,T]}(\phi) ~=~ -\inf_{\phi\in F:\phi(0)=z} I_{[0,T]}^+(\phi).
\end{align*}
This proves that the sequence of rescaled processes $(Z_+^\eps(t))$ satisfies SPLD upper
bound with the rate function $I_{[0,T]}^+$.

To get the SPLD lower bound it is sufficient to show  that 
\be\label{e4-5}
\lim_{\delta'\to 0} ~\lim_{\delta\to 0} \;\liminf_{\eps\to 0} \; \inf_{z: |z-\phi(0)|<\delta}
\eps \log\P_{[z/\eps]}\left( \sup_{t\in[0,T]} |Z^\eps_+(t) - \phi(t)| < \delta' \right) 
\geq~ - I_{[0,T]}(\phi) 
\ee 
for every  absolutely continuous function $\phi :~[0,T]\to\R^{d-1}\times\R_+$. If $\phi(t)$ belongs
to the interior of the half-space $\R^{d-1}\times\R_+$ for every $t\in[0,T]$, then 
this inequality follows from Mogulskii's theorem because in this case, for $\delta' >0$ 
small enough, the process $(Z(t))$ does 
    not exit from $\Z^{d-1}\times\N^*$ before the 
    time $T/\eps$ and $Z^\eps_+(t)=Z^\eps(t)$ for all $t \leq T$  whenever 
\[
\sup_{t\in[0,T]} |Z^\eps(t) - \phi(t)| < \delta'.
\]
Let now $\phi :~[0,T]\to\R^{d-1}\times\R_+$ be an arbitrary absolutely continuous function. For
$\sigma >0$, consider the function $\phi_\sigma$ on $[0,T]$ defined by 
\[
\phi_\sigma(t) ~=~ \phi(t) ~+~ \sigma e_d
\]
where $e_d = (0,\ldots,0,1)\in\Z^d$. Because of Lemma~\ref{lem1-4}, there are $C\geq 1$ and
$\theta >0$ such that for any $z\in\R^d$
satisfying the inequality $|\phi(0) - z| < \delta$ and for any $\eps >0$ there is a
sequence of points $z_0, z_1, \ldots,z_N\in\Z^{d-1}\times\N^*$ with $z_0 = [z/\eps]$ and
$z_N=[\phi_\sigma(0)/\eps]$ for which  the following inequalities  hold
\begin{align}
N &~\leq~ C\left| [z/\eps] - [\phi_\sigma(0)/\eps]\right| ~\leq~ 2 d C + |z -
\phi_\sigma(0)|/\eps \nonumber \\ 
&~\leq~ 2 d C + (\delta + \sigma)/\eps, \label{e4-6}
\end{align}
\be\label{e4-7}
 |z_i-z_{i-1}| ~\leq~ C \quad \text{and} \quad  p(z_{i-1},z_i) ~\geq~ \theta, \quad \forall \; i=1,\ldots, N. 
\ee 
Hence, on the event 
\[
E_{\sigma,\sigma',\eps} ~=~ \left\{Z_+(0) = z_0, \ldots, Z_+(N)=z_N\right\} \cap\left\{ \sup_{t\in[0,T]} \left|
Z_+^\eps(t+N\eps) - \phi_\sigma(t)\right| < \sigma' \right\}, 
\]
for $t\in[0,N\eps]$, using the inequality \eqref{e4-6} and the first inequalit of \eqref{e4-7} we obtain 
\begin{align*}
\left|Z_+^\eps(t) - \phi(t)\right| &~\leq~ \left| Z_+^\eps(t) - \phi_\sigma(0)\right|
+ |\phi_\sigma(0) - \phi(0)| + |\phi(0)-\phi(t)| \\ 
&~\leq~ \left| Z_+^\eps(t) - \phi_\sigma(0)\right|
+ \sigma + \max_{t\in[0, N\eps ]} |\phi(0)-\phi(t)| \\ 
&~\leq~ \max_{i=1,\ldots,N} \left|
\eps z_i - \phi_\sigma(0)\right| + \sigma + \max_{t\in[0, N\eps ]} |\phi(0)-\phi(t)| \\ 
&~\leq~ \eps \max_{i=1,\ldots,N} \left|
z_i - z_N\right| +  \eps d + \sigma + \max_{t\in[0, N\eps ]} |\phi(0)-\phi(t)| \\  
&~\leq~ \eps(N C + d) + \sigma + \max_{t\in[0, N\eps ]} |\phi(0)-\phi(t)| \\
&~\leq~ 2 d C^2\eps + (\delta + \sigma) C + \eps d + \sigma + \max_{\substack{s,s'\in[0,
      T] : \\|s'-s| \leq  \delta + \sigma + 2 d C\eps}} 
|\phi(s)-\phi(s')|  
\end{align*}
and for $t\in[N\eps, T]$ we get 
\begin{align*}
\left|Z_+^\eps(t) - \phi(t)\right| &\leq  \left|Z_+^\eps(t) - \phi_\sigma(t-
N\eps)\right|  + |\phi_\sigma(t- N\eps)  - \phi(t)|\\ 
&\leq  \sigma' + |\phi_\sigma(t- N\eps) -  \phi(t)|\\ 
&\leq  \sigma' + \sigma + |\phi(t-N\eps) - \phi(t)|\\ 
&\leq \sigma' + \sigma + \max_{\substack{s,s'\in[0, T] : \\|s'-s| \leq  \delta + \sigma + 2 d C\eps}}
|\phi(s)-\phi(s')| 
\end{align*} 
where 
\[
\max_{\substack{s,s'\in[0, T] : \\|s'-s| \leq  \delta + \sigma + 2 d C\eps}}
|\phi(s)-\phi(s')| ~\to~ 0 
\]
when $\delta, \sigma, \eps \to 0$ because the function $\phi$ is continuous on
$[0,T]$. These inequalities show that for any $\delta' >0$ there is $\sigma_0 >0$ such that
for all $\sigma, \delta, \sigma', \eps \in ]0,\sigma_0[$, on the event
    $E_{\sigma,\sigma',\eps}$ the following inequality holds 
\[
\sup_{t\in[0,T]} \left|
Z_+^\eps(t) - \phi(t)\right| < \delta' 
\]
Using the second inequality of \eqref{e4-7} and Markov property  from this it follows that
for all $0 < \eps,\delta,
\sigma' < \sigma < \sigma_0$,  
\begin{align*}
\eps \log \P_{[z/\eps]}&\left( \sup_{t\in[0,T]} |Z^\eps_+(t) - \phi(t)| < \delta' \right)
\\ &~\geq~ N\eps \log \theta + \eps \log ~\P_{[\phi_\sigma(0)/\eps]}\left( \sup_{t\in[0,T]} |Z^\eps_+(t) -
\phi_\sigma(t)| < \sigma' \right)\\
&~\geq~ (2d C \eps + \delta + \sigma) \log \theta + \eps \log
~\P_{[\phi_\sigma(0)/\eps]}\left( \sup_{t\in[0,T]} |Z^\eps_+(t) - 
\phi_\sigma(t)| < \sigma' \right) 
\end{align*}
and consequently, 
\begin{align}
\liminf_{\eps\to 0} &~\eps\log~\P_{[z/\eps]}\left( \sup_{t\in[0,T]} |Z^\eps_+(t) - \phi(t)| <
\delta' \right) \nonumber\\
&\geq~  (\delta +\sigma) \log \theta ~+~ \liminf_{\eps\to 0} \eps~\log~ \P_{[\phi_\sigma(0)/\eps]}\left(
\sup_{t\in[0,T]} |Z^\eps_+(t) - 
\phi_\sigma(t)| < \sigma' \right) \nonumber\\
&\geq~  (\delta + \sigma) \log \theta ~+~ \liminf_{\eps\to 0} \eps~\log~ \P_{[\phi_\sigma(0)/\eps]}\left(
\sup_{t\in[0,T]} |Z^\eps(t) - 
\phi_\sigma(t)| < \sigma' \right) \label{e4-8}
\end{align}
where the last inequality holds because for $0<\sigma' < \sigma$, on the event 
\[
\left\{ \sup_{t\in[0,T]} |Z^\eps(t) - 
\phi_\sigma(t)| < \sigma'\right\},
\]
we have the equality $Z_+(t) = Z(t)$ for all $t\in[0,T]$. Relation \eqref{e4-8} combined
with SPLD lower bound 
of Mogulskii's theorem shows that the left hand side of \eqref{e4-5} is greater than 
\begin{align*}
\limsup_{\sigma\to 0} \!\!\lim_{\sigma'\to 0} \liminf_{\eps\to 0} \eps \log \P_{[\phi_\sigma(0)/\eps]}&\!\left(
\sup_{t\in[0,T]} |Z^\eps(t) - 
\phi_\sigma(t)| < \sigma' \right) \!\geq - \liminf_{\sigma\to 0} I_{[0,T]}(\phi_\sigma)
\end{align*}
and hence, using the equalities 
\[
I_{[0,T]}(\phi_\sigma) ~=~ \int_0^T (\log \varphi)^*(\dot{\phi}_\sigma(t)) \, dt ~=~ \int_0^T
(\log \varphi)^*(\dot{\phi}(t)) \, dt ~=~ I_{[0,T]}(\phi)  
\]
we obtain \eqref{e4-5}
\end{proof}

\subsection{Large deviation estimates for the Green's functions} The lower rough
logarithmic  estimates for the Green's function $G_+(z,z')$ are derived now from the SPLD properties of
the scaled processes 
$Z^\eps_+(t)=\eps Z_+([t/\eps])$.

\begin{prop}\label{pr4-2} Under the hypotheses (A), for any $q\in {\cal S}_+^d$ and 
  $z\in\Z^{d-1}\times\N^*$,  
\be\label{e4-9}
\liminf_{n\to\infty} ~\frac{1}{|z_n|}\log G_+(z,z_n) ~\geq~ -  a(q)\cdot q 
\ee
when  $|z_n|\to\infty$ and $z_n/|z_n|\to q$ as $n\to\infty$, with $z_n\in\Z^{d-1}\times\N^*$.
\end{prop}
 \begin{proof} 
Indeed,  if the hypotheses (A) are satisfied
  then by Lemma~\ref{lem1-4}, the Markov process $(Z_+(t))$ satisfies communication
  condition (A1') and hence,  for any $z,z'\in \Z^{d-1}\times\N^*$, if $z\not= z'$ then
  there is $0<t\leq C|z-z'|$ such that 
\[
\P_z(Z_+(t) = z') ~\geq~ \theta^t ~\geq~ \theta^{C|z'-z|}.
\]
Using the inequality  $G_+(z,z_n) \geq G_+(z,z') \,\P_{z'}(Z_+(t) = z_n) $ for
  $z'\not= z_n$, we get therefore 
\[
G_+\bigl(z,z_n\bigr) ~\geq~   G_+\bigl(z,z'\bigr)\theta^{C|z'-z_n|},  \quad \quad \quad \forall
n\in\N, \; z'\in \Z^{d-1}\times\N^*. 
\]
Using moreover the inequality 
\[
\text{Card}\{z\in\Z^d : |z - a|<\delta R\} ~\leq~ (2\delta R + 1)^{d}
\]
with $R = |z_n|$ and $a = |z_n| q$ we obtain 
\[
G_+\bigl(z,z_n\bigr) ~\geq~ \frac{1}{(2\delta |z_n|+1)^{d}} \theta^{2C\delta |z_n|}
\sum_{\substack{z' :~|q - z'/|z_n|| <\delta}}
G_+\bigl(z,z'\bigr) 
\]
for all those $n\in\N$ for which  $|q
- z_n/|z_n|| <\delta$ and consequently, 
\begin{align*}
\lim_{n\to\infty} ~\frac{1}{|z_n|}\log G_+(z,z_n)  &~\geq~ - 2C \delta \log \theta +
\liminf_{n\to\infty} ~\frac{1}{|z_n|} 
~\log \sum_{\substack{z' :~|q - z'/|z_n||<\delta }} G_+\bigl(z,z'\bigr)\\
&~\geq~ - 2C \delta \log \theta +
\liminf_{\eps\to 0} ~\eps 
~\log \sum_{\substack{z' :~|q - \eps z'|<\delta }} 
G_+\bigl(z,z'\bigr). 
\end{align*}
Finally, letting $\delta\to 0$ and using SPLD lower bound for the family of scaled process
$Z^\eps_+(t) = \eps Z_+([t/\eps])$  we conclude that  
\begin{align*}
\lim_{n\to\infty} ~\frac{1}{|z_n|}\log G_+(z,z_n)  &~\geq~ 
\lim_{\delta\to 0} \liminf_{\eps\to 0} ~\eps ~\log \sum_{\substack{z' : ~|q - \eps z'|<\delta }}
G_+\bigl(z,z'\bigr) \nonumber\\ &~\geq~  \lim_{\delta\to 0}  \liminf_{\eps\to 0} \eps \log
\P_{z}\Bigl(  |Z_+^\eps(T) - q| < 
\delta \Bigr) \nonumber 
\\  &~\geq~ - \inf_{\phi :~ \phi(0)=0,\, \phi(T)=q}
I_{[0,T]}^+(\phi)
\end{align*}
for every $T>0$. To get \eqref{e4-9}  it is sufficient now to notice  that the right hand side of the last
inequality is greater than $-a(q)\cdot 
q$ because for $T=1/|\nabla\varphi(a(q))|$ and for the linear function $
\phi(t) = vt$ with $v = \nabla\varphi(a(q)) = q |\nabla\varphi(a(q))|$, from \eqref{e4-4}
it follows that  
\[
I_{[0,T]}^+(\phi) ~=~ T (\log \varphi)^*(v) ~=~  T a(q)\cdot v ~=~ a(q)\cdot
q. 
\]
\end{proof}

\section{Harmonic functions}  
The harmonic functions of the Markov process $(Z_+(t))$ are now identified.
The main result of this section is the following proposition.  

\begin{prop}\label{pr5-1} Under the hypotheses (A), the following assertions hold. 

\noindent 
1) A non-negative  function $h$ is harmonic for the Markov process $(Z_+(t))$ if and only if
 there  is a positive measure $\nu_h$ on $\partial_+ D$ such that 
\be\label{e5-1}
h(z) ~=~ \int_{\partial_+ D} h_{a,+}(z) \, d\nu_h(a), \quad \quad \forall
z\in\N^*\times\Z^{d-1}.  
\ee
2) For every $a=(\a,\beta)\in\partial_+ D$ with $\a\in\R^{d-1}$ and $\beta\in\R$, the
constant multiples of of the function $h_{a,+}$ are the 
only non-negative harmonic functions for which 
\be\label{e5-2}
\sup_{x\in\R^{d-1}} \exp(- \a\cdot x)  h(x,y) ~<~ +\infty, \quad \quad \forall
y\in\N^*. 
\ee
3) The  constant multiples of the functions $h_{a,+}$ with
  $a\in\partial_+ D$, are the only minimal harmonic  
functions of the Markov process $(Z_+(t))$. 
\end{prop}
Throughout
this section  the following notations are used : for 
$a=(\a,\beta)\in\R^d$ we denote by $\a$  [resp. $\beta$] the vector of $d-1$ first
coordinates  [resp. last coordinate] of $a$. For $z=(x,y)\in\Z^d$ the variavles $x$ and $y$ are
defined in the similar way. 

To prove Proposition~\ref{pr5-1} we
combine the properties of Markov-additive processes and the results of
Doney~\cite{Doney:02}. The Markov process $(Z_+(t)) = (X_+(t),Y_+(t))$ is
Markov-additive  with an 
additive part $X_+(t)$ taking the values in $Z^{d-1}$ and a Markovian part $Y_+(t)$ taking
the values in $\N^*$.   Under the hypotheses~$(A)$, the Markovian part
$Y_+(t)$ is irreducible on $\N^*$ and hence, the Feynman-Kac transfom matrix ${\cal
  P}_+(\alpha) = \bigl( {\cal P}_+(\alpha, 
y,y'), \; y,y'\in\N^*\bigr)$ with 
\[
{\cal P}_+(\alpha, y,y') ~=~ \E_{(0,y)}\bigl( \exp(\alpha\cdot X_+(1)) ; \, Y_+(1) = y'\bigr)
\]
is also irreducible. The quantity
  $e^{\lambda_+(\alpha)}$ with 
\[
\lambda_+(\alpha) ~\dot=~ \limsup_n \frac{1}{n} \log {\cal P}^{(n)}_+(\alpha, y,y') 
\]
is the {\it spectral radius} of the
  transform matrix ${\cal P}_+(\alpha)$. Recall that this limit does not depend on
  $y,y'\in \N^*$ because the matrix ${\cal P}_+(\alpha)$ is irreducible
  (see~\cite{Seneta}). 

To identify the harmonic functions of the Markov process $(Z_+(t))$ we identify first the
function $\lambda_+(\cdot)$. This is a subject of the following lemma. 

\begin{lemma}\label{lem5-1} Under the hypotheses $(A)$, 
\[
\lambda_+(\alpha) ~=~ \inf_{\beta\in\R} \log \varphi(\alpha,\beta), \quad \quad \forall
\a\in\R^{d-1}.
\]
\end{lemma}
\begin{proof} To prove this lemma we consider  a random walk $(Z(t)) = (X(t),Y(t))$ on
  $\Z^d$ having transition probabilities $p(z,z')=\mu(z'-z)$ as a Markov-additive process
  on $\Z^{d-1}\times\Z$ with an additive part $(X(t))$ on $\Z^{d-1}$ and a Markovian part
  $(Y(t))$ on $Z$. For such a Markov-additive process, the Feynman-Kac transfom matrix ${\cal
  P}(\alpha) = \bigl( {\cal P}(\alpha, 
y,y'), \; y,y'\in\Z\bigr)$ is defined by 
\[
{\cal P}(\alpha, y,y') ~=~ \E_{(0,y)}\bigl( \exp(\alpha\cdot X(1)) ; \; Y(1) = y'\bigr), \quad
\quad \forall y,y\in\Z
\] 
and its spectral radius $e^{\lambda(\a)}$ is given by 
\[
\lambda(\alpha) ~\dot=~ \limsup_{t\to\infty} \frac{1}{t} \log \P_{(0,y)}\bigl( \exp(\alpha\cdot X(t)), \;
Y(t) = y'\bigr).
\]
The first step of our proof shows that 
\be\label{e5-3}
\lambda(\alpha) ~=~ \inf_{\beta\in\R} \log \varphi(\alpha,\beta).
\ee
Indeed, under the hypotheses $(A)$, the function $\varphi$ is
  convex and has compact level sets. For every $\alpha$, there
  is therefore  $\beta_\alpha^0\in\R$ such that 
\[
\inf_{\beta\in\R} \varphi(\alpha,\beta) = \varphi(\alpha,\beta_\alpha^0).
\]
A twisted random walk $(\tilde{Y}(t))$ on $\Z$ with transition
probabilities 
\[
\tilde{p}\bigl(y,y'\bigr) ~=~ \sum_{x\in\Z^{d-1}} \mu(x,y'-y) \exp\bigl(\alpha\cdot x +
    \beta_\alpha^0 (y'-y)\bigr)/ \varphi(\alpha,\beta_\alpha^0) 
\]
has finite variance and zero mean because 
\[
\E(\tilde{Y}(1)) ~=~ \frac{\partial}{\partial\beta}\log \varphi(\alpha,\beta_\alpha^0) ~=~ 0,
\]
from which it follows that  
\[
\limsup_{t\to\infty} \frac{1}{t} \log \tilde{p}^{(t)}(y,y') ~=~ 0
\]
for all $y,y'\in\Z$ (see Spitzer~\cite{Spitzer}). The last relation combined with the
equality 
\[
\E_{(0,y)}\Bigl(\exp\bigl(\alpha\cdot X(t)\bigr), \; Y(t) =
    y'\Bigr)  ~=~ \left(\varphi(\alpha,\beta_\alpha^0)\right)^t \exp\bigl(\beta_\alpha^0(y-y')\bigr) ~\tilde{p}^{(t)}(y,y). 
\]
proves  \eqref{e5-3}. 

Now, to complete the proof of Lemma~\ref{lem5-1}, we show that $
\lambda_+(\alpha)=\lambda(\alpha)$. 
For this we notice that  ${\cal P}(\alpha; y,y') = {\cal P}_+(\alpha; y,y')$ for all
$y,y'\in\N^*$. Since under the Assumption $(A)$, the matrices ${\cal
  P}_+(\alpha) = ({\cal P}(\a; y,y'), \; y,y'\in\N^*)$ and ${\cal
  P}(\alpha) = ({\cal P}(\a; y,y'), \; y,y'\in\Z)$ are irreducible then by Theorem~6.3 of
 Seneta~\cite{Seneta} (see also Proposition~2 of Ignatiouk~\cite{Ignatiouk:02}), 
\[
\lambda(\alpha) ~=~ \sup_{K\subset\subset \Z} ~\lambda_K(\alpha) \quad \text{and} \quad
\lambda_+(\alpha) ~=~ \sup_{K\subset\subset \N^*} ~\lambda_K(\alpha) 
\]
where  both supremums are taken over all finite subsets $K$ and for any finite set $K\subset \Z$,
  $\exp(\lambda_K(\alpha))$ is the maximal real eigenvalue of the truncated matrix $\bigl({\cal
  P}(\alpha; y,y'); \; y,y'\in K\bigr)$. Moreover, since ${\cal
  P}(\alpha; y'+y,y''+y) = {\cal P}(\alpha; y',y'')$ for all $y,y',y''\in\Z$, then for
  every finite set $K\subset \Z$ we have also 
\[
\lambda_K(\alpha) ~=~ \lambda_{K+y}(\alpha), \quad \quad \forall y\in\Z
\] 
and consequently, $\lambda_+(\alpha)=\lambda(\alpha)$. 
\end{proof}

\bigskip

Remark that by Lemma~\ref{lem5-1}, $\lambda_+(\a) \leq 0$ if and only if
$\varphi(\alpha,\beta) \leq 1$ for some $\beta\in\R$ and hence, 
the mapping $a=(\a,\beta) \to \a$ determines a one to
one and on correspondence from $\partial_+ D$ to $
\{ \a\in\R^{d-1} :~ \lambda_+(\a) \leq 0\}$. 
Using Lemma~\ref{lem5-1} and Proposition~\ref{pr3-0} we obtain therefore the following
statement. 

\begin{lemma}\label{lem5-2}  Under the Assumptions $(A)$,   every minimal harmonic function
  $h$ of the Markov process $(Z_+(t))$ is of the form   
\be\label{e5-5}
h(x,y) = \exp\bigl(\alpha \cdot x\bigr) h(0,y) ~>~ 0 , \quad \quad \forall (x,y)\in\Z^{d-1}\times\N^*
\ee
for some $a=(\a,\beta)\in\partial_+ D$. 
\end{lemma}

To identify the harmonic functions 
satisfying the equation  \eqref{e5-5}  we need the following lemma. 

\begin{lemma}\label{lem5-3} For every aperiodic random walk $(Y(t))$ on $\Z$ having transition
  probabilities $P(y,y')=P(0,y'-y)$ such that for some $\delta < 0$, 
\[
\sum_y e^{-\delta y} P(0,y) ~<~ \infty, \quad \quad \sum_y |y| P(0,y) ~<~ \infty, \quad \text{and} \quad m
~\dot=~ \sum_y yP(0,y) \geq 0,   
\]
the constant multiples of the function 
\be\label{e5-6}
f(y) ~=~ \begin{cases}  \P_{y}( \tau = +\infty) &\text{if
   $ m > 0$,}\\
y - \E_y(Y(\tau)) &\text{if
   $ m = 0,$}
\end{cases}
\ee
with $\tau = \inf\{t\geq 0 : Y(t)\leq 0\}$, are the only positive solutions of the
equation 
\be\label{e5-7}
\sum_{y'>0} P(y,y') f(y') ~=~ f(y), \quad y\in\N^*.
\ee
\end{lemma}
\begin{proof} This statement follows from Theorem~1 of Doney~\cite{Doney:02} (when $m=0$
  this is a consequence of Example E~27.3 in Chapter VI of
  Spitzer~\cite{Spitzer}). This theorem proves that for an aperiodic random
  walk $(Y(t))$ on $\Z$ having  transition
  probabilities $P(y,y')=P(0,y'-y)$ and a non-negative mean $m\geq 0$, 
the only positive solutions of the
equation 
\[
\sum_{y'\geq 0} P(y,y') f(y') ~=~ f(y), \quad y\in\N.
\]
 are the constant multiples of the renewal function of strict increasing ladder heights of
 $(-Y(t))$. Hence, to prove Lemma~\ref{lem5-3} it is sufficient to show that the function
 \eqref{e5-6} is well defined and satisfies the equation \eqref{e5-7}. 

When $m\not=0$ the function $f$ is
  clearly positive and well defined. Suppose now that $m=0$ and let us consider the
  function 
\[
g_+(y,y') ~=~ \sum_{t=0}^\infty \P_y(Y(t) = y', \; \tau > t).
\]
 In this case, for any $\delta >0$,
\begin{align*}
0 ~\leq~ - \E_y(Y(\tau)) &~\leq~ \frac{1}{\delta} ~\E_y\left(e^{-\delta Y(\tau)}\right) ~=~
  \sum_{y' > 0} g_+(y,y') \sum_{y'' \leq 0} e^{-\delta y''}P(y',y'')\\ 
&~\leq~   \sum_{y''} e^{-\delta y''}P(0,y'') ~\sum_{y' > 0} g_+(y,y') e^{-\delta y'}\\ 
&~\leq~   \sum_{y''} e^{-\delta y''}P(0,y'') ~g_+(y,y)  ~\sum_{y' > 0} e^{-\delta y'} 
\end{align*}
where the last inequality holds because  for the ``reversed'' random walk $\tilde{Y}(t)$
with transition probabilities 
$\tilde{P}(y,y') = P(0, y-y') = P(y',y)$, 
\begin{align*}
g_+(y,y') &~=~ \sum_{t=0}^\infty \P_y(Y(t) = y', \; \tau > t) ~=~ \sum_{t=0}^\infty
\P_{y'}(\tilde{Y}(t) = y, \; \tau > t) \\ &~\leq~ \P_{y'}( \tilde{Y}(t) = y \; \text{ for
  some } \; t \geq 0) 
~\sum_{t=0}^\infty \P_{y}(\tilde{Y}(t) = y, \; \tau > t) 
\\ &~\leq~  
~\sum_{t=0}^\infty \P_{y}(\tilde{Y}(t) = y, \; \tau > t) ~=~ g_+(y,y) 
\end{align*}
for every $y'\in\N^*$. This proves that the function $f$ is positive and finite. 

Finally, straightforward calculations show that the function $f$ 
  satisfies the equation \eqref{e5-7}. \end{proof}

\medskip

The harmonic functions 
satisfying the equation  \eqref{e5-5} are now identified.  
\begin{lemma}\label{lem5-4} Under the Assumptions $(A)$, for every $a=(\alpha,\beta)\in\partial_+ D$, 
 the constant multiples of $h_{a,+}$  are the only harmonic 
functions of the Markov process $(Z_+(t))$ for which  \eqref{e5-5} holds with a given $\a$.
\end{lemma}
\begin{proof} To prove this lemma it is sufficient to show 
that  for every $a\in\partial_+
D$, the constant multiples of $\tilde{h}_{a,+}(z) = \exp(-a\cdot z) h_{a,+}(z)$ are the only
positive solutions of the equation 
\be\label{e5-8}
\sum_{z'\in\Z^{d-1}\times\N^*} p(z,z') \exp\bigl(a\cdot (z'-z)\bigr) \tilde{h}(z') ~=~
\tilde{h}(z), \quad \quad z\in\Z^{d-1}\times\N^*
\ee
satisfying the equality $\tilde{h}(x,y) ~=~ \tilde{h}(0,y)$ for all
$(x,y)\in\Z^{d-1}\times\N^*$.
 To prove such a property 
for a given $a=(\a,\beta)\in\partial_+ D$, we consider a twisted random walk
$\tilde{Z}(t)=(\tilde{X}(t),\tilde{Y}(t))$ on $\Z^d$ with transition probabilities
\[
\tilde{p}(z,z') ~=~ \exp\bigl(a \cdot(z'-z)\bigr) \mu(z'-z).
\]
Then  for every $z=(x,y)\in\Z^{d-1}\times\N^*$,  
\[
\tilde{h}_{a,+}(z) ~=~ f_a(y) ~\dot=~ \begin{cases}  \P_{y}( \tau_a = +\infty) &\text{if
   $a\in\partial_+ D\setminus\partial_0 D$}\\
y - \E_y(\tilde{Y}(\tau_a)) &\text{if
   $ a\in\partial_0 D$}
\end{cases}
\]
where  $\tau_a = \inf\{t\geq 0 : \tilde{Y}(t)
\leq 0\}$. Moreover, because of the assumption (A), the last coordinate 
$(\tilde{Y}(t))$ of $(\tilde{Z}(t))$ is an aperiodic random walk on $\Z$ with transition probabilities 
\[
P_a(y,y') ~=~ P_a(0,y'-y) ~=~ \sum_{x\in\Z^{d-1}}
  \mu(x,y'-y) \exp\bigl(\alpha\cdot x + \beta(y'-y)\bigr)
\]
and mean
\[
m_a ~\dot=~ \E_0(\tilde{Y}(1)) ~=~
  \sum_{z=(x,y)\in\Z^d} y \,\mu(z) 
  \exp(a\cdot z) ~=~ \frac{\partial}{\partial\beta}
  \varphi(\alpha,\beta) ~\geq~ 0  
\]
where $m_a = 0$ if and only if $a\in\partial_0 D$. Hence, by Lemma~\ref{lem5-3}, the
constant multiples of the function
$f_a$  are  the only positive solution of the equation  
\[
\sum_{y' > 0} P_a(y,y') f(y') = f(y), \quad \quad  y\in\N^* 
\] 
and therefore, the
constant multiples of $\tilde{h}_{a,+}$ are the only positive 
solutions of the equation \eqref{e5-8} satisfying the equality $\tilde{h}(x,y) = \tilde{h}(0,y)$ for all
$(x,y)\in\Z^{d-1}\times\N^*$. 
\end {proof}

\medskip

The last lemma combined with Lemma~\ref{lem5-2} proves the following statement. 

\begin{lemma}\label{lem5-5}  Under the hypotheses $(A)$,   every minimal harmonic function
  of the Markov process $(Z_+(t))$ is of the form   $h = ch_{a,+}$ 
with some $c>0$ and $a\in\partial_+ D$. 
\end{lemma}

\noindent
{\bf Proof of Proposition~\ref{pr5-1}}. We are ready now to prove the 
 representation \eqref{e5-1}.  By the Poisson-Martin representation theorem
(see Woess~\cite{Woess}), every 
non-negative harmonic function of the Markov process $(Z_+(t))$ is of the form 
\be\label{e5-9}
h(z) ~=~ \int_{\partial_m(\Z^{d-1}\times\N^*)} K_+(z,\gamma) \, d\tilde\nu_h(\gamma), \quad \quad
\forall z\in\Z^{d-1}\times\N^*
\ee
with some Borel measure $\tilde\nu_h\geq 0$ on the minimal Martin boundary
$\partial_m(\Z^{d-1}\times\N^*)$.   Recall that 
 $K_+(z,\gamma)$ is
the Martin kernel of the Markov process $(Z_+(t))$. The Martin
compactification of $\Z^{d-1}\times\N^*$ for the Markov process $(Z_+(t))$ is the unique
smallest compactification of 
the  set $\Z^{d-1}\times\N^*$ for  which the  functions 
\[
z' \to K_+(z,z') ~=~ G_+(z,z')/G_+(z_0,z')
\]
with a given $z_0\in\Z^{d-1}\times\N^*$ extend continuously for all
$z\in\Z^{d-1}\times\N^*$. 
The mapping $\gamma\to
K_+(z,\gamma)$ is therefore continuous on $\partial_m(\Z^{d-1}\times\N^*)$ for every
$z\in\Z^{d-1}\times\N^*$. Moreover, for
every $\gamma\in\partial_m(\Z^{d-1}\times\N^*)$, according to the definition of the minimal Martin boundary (see Woess~\cite{Woess}), the function $z \to K(z,\gamma)$
is a minimal harmonic function for the Markov process $(Z_+(t))$ with $K(z_0,\gamma) = 1$. By
Lemma~\ref{lem5-5}, from this it follows that 
\be\label{e5-10}
K_+(z,\gamma) = c_\gamma h_{a(\gamma),+}(z),\quad \quad \forall z\in\Z^{d-1}\times\N^*
\ee
with some $a(\gamma)\in\partial_+ D$ and $c_\gamma= 1/h_{a(\gamma),+}(z_0)$. Hence, for
$z_0=(x_0,y_0)$ and $z=(x,y_0)$,  
\[
K_+(z,\gamma) ~=~ \exp\bigl(\a(\gamma)\cdot(x-x_0)\bigr)
\]
where $\a(\gamma)\in\R^{d-1}$ denotes the vector of $d-1$ first coordinates of
$a(\gamma)$. Since the mapping $\gamma\to
K_+(z,\gamma)$ is continuous on $\partial_m(\Z^{d-1}\times\N^*)$ for every
$\Z^{d-1}\times\N^*$, then  from the last relation it follows  that the mapping
$\gamma\to\a(\gamma)\cdot x$ is continuous on 
$\partial_m(\Z^{d-1}\times\N^*)$ for
every $x\in\Z^{d-1}$ and consequently, the mapping $\gamma\to\a(\gamma)$ is also continuous on
$\partial_m(\Z^{d-1}\times\N^*)$.  From this it follows that the mapping
$\gamma\to a(\gamma)$ from $\partial_m(\Z^{d-1}\times\N^*)$ to $\partial_+ D$ is also
continuous  because for every $\a\in\R^{d-1}$ for which $\inf_{\beta\in\R}
\varphi(\a,\beta)\leq 1$  there is a unique 
$\beta_\a\in\R$ such that $(\a,\beta_a)\in\partial_+ D$, and the mapping $\a\to
(\a,\beta_\a)$ from 
$\{\a\in\R^{d-1} : \inf_\beta
\varphi(\a,\beta)\leq 1\}$ to $\partial_+ D$ is continuous. The equalities 
\eqref{e5-9} and \eqref{e5-10} imply therefore \eqref{e5-1} with the positive Borel measure $\nu_h$ on
$\partial_+ D$ defined by
\[
\nu_h(B) ~=~ \int_{\{\gamma : a(\gamma)\in B\}} c_\gamma \, d\tilde\nu_h(\gamma) 
\]
for every Borel subset $B\subset\partial_+ D$ (the set $\{\gamma : a(\gamma)\in B\}$ is
here measurable because the mapping $\gamma \to a(\gamma)$ is continuous and hence, the
measure $\nu$ is well defined).  The first assertion of Proposition~\ref{pr5-1} is therefore proved.

Now, to prove the second assertion it is sufficient to show that a non-zero harmonic
function $h\geq 0$ satisfies  \eqref{e5-2} for some $a=(\a,\beta)\in\partial_+ D$ if and
only if $$\supp(\nu_h) =\{a\}.$$ By Lemma~\ref{lem5-4}, for every $a\in\partial_+ D$, the function $h_{a,+}$
is harmonic for the Markov process $(Z_+(t))$ and satisfies  \eqref{e5-2} because for every
$z=(x,y)\in\Z^{d-1}\times\N^*$, 
\[
\exp(-\a\cdot x) h_{a,+}(x,y) ~=~ h_{a,+}(0,y). 
\] 
Conversely, if  $\supp(\nu_h) \not=\{\hat{a}\}$ for some $\hat{a}=(\hat\a,\hat\beta)\in\partial_+ D$,
then there are $a_0=(\a_0,\beta_0)\in \partial_+ D$ and $\eps >0$ such that $|\hat{a} - a_0|
\geq |\hat{\a} - \a_0| >
\eps$ and $$\nu(B(a_0,\eps)\cap\partial_+ D) > 0$$ where $B(a_0,\eps)$ denotes the open
ball in $\R^d$ 
centered at $a_0$ and having the radius $\eps>0$.  Moreover, since for every
$a=(\a,\beta)\in B(a_0,\eps)$,
\[
|\a - \a_0| ~\leq~ |a - a_0| ~<~ \eps ~<~ |\hat{\a} - \a_0|
\]
then there exists 
 $x_0\in\Z^{d-1}$   such that $(\a
- \hat\a)\cdot x_0 > 0$ for all $a=(\a,\beta)\in B(a_0,\eps)$ and consequently, using
 Fatou lemma we get 
\begin{align*}
\sup_{x\in\R^{d-1}} e^{-\hat{\a}\cdot x} h(x,y)  &~\geq~ \limsup_{n\to\infty} ~e^{- n\hat{\a}\cdot
  x_0} h(n x_0,y) \\ &~\geq~   \limsup_{n\to\infty} ~\int_{B(a_0,\eps)\cap\partial_+
  D} e^{n (\a - \hat\a)\cdot x_0} h_{a,+}(0,y)  \, d\nu_h(a) \\ &~\geq~
  ~\int_{B(a_0,\eps)\cap\partial_+ 
  D} \lim_{n\to\infty} e^{n (\a - \hat\a)\cdot x_0} h_{a,+}(0,y)  \, d\nu_h(a)   ~=~
 + \infty. 
\end{align*}
Hence, if a non-zero harmonic
function $h\geq 0$ satisfies  \eqref{e5-2} for some $a=(\a,\beta)\in\partial_+ D$  then 
 $supp(\nu ) = \{a\}$ and consequently $h
= c h_{a,+}$ with some $c >0$. The second assertion of Proposition~\ref{pr5-1} is proved. 

Finally,  by Lemma~\ref{lem5-4}, for every
$a\in\partial_+D$, the function $h_{a,+} > 0$ is 
harmonic for the Markov process $(Z_+(t))$. Moreover, if a non-negative harmonic function $h$ satisfies 
the inequality $h \leq h_{a,+}$ for some $a\in\partial_+ D$ then for $h$ the inequality
\eqref{e5-2} holds with the same $a\in\partial_+ D$ and
consequently $h=ch_{a,+}$ for some $c\geq 0$. For every $a\in\partial_+D$, the harmonic
function $h_{a,+} > 0$ is therefore minimal and  conversely, by Lemma~\ref{lem5-5}, every
minimal harmonic function of the Markov process $(Z_+(t))$ is of the form $c h_{a,+}$ with
some $c>0$ and $a\in\partial_+ D$. Proposition~\ref{pr5-1} is proved.

\section{Proof of Theorem~\ref{th1-4}} The first assertion of Theorem~\ref{th1-4} is
proved by Proposition~\ref{pr5-1}. To prove the second assertion of
Theorem~\ref{th1-4} we use
Proposition~\ref{pr3-2} combined with Proposition~\ref{pr4-2} and
Proposition~\ref{pr5-1}.  

Consider a sequence of points $z_n \in\Z^{d-1}\times\N^*$ with $\lim_n |z_n| = \infty$
such that $\lim_n z_n/|z_n| = q\in{\cal S}_+^d$.  Suppose first that $m/|m|\in{\cal
  S}^d_+$ and let $q = m/|m|$. Then $q=q(0)$, $a(q)=0$ and 
by Proposition~\ref{pr4-2}, 
\be\label{4-10}
\liminf_{n\to\infty} \frac{1}{|z_n|} \log  G_+(z,z_n) ~\geq~
- a(q)\cdot q = 0.  
\ee 
Furthermore, by Lemma~\ref{lem1-4}, the Markov process $(Z_+(t))$ satisfies the
communication condition (A1'). The function 
\[
\varphi(a) ~=~ \sup_{z\in\Z^{d-1}\times\N^*} \E_z\left(\exp\bigl(a\cdot (Z_+(1)-z)\bigr)\right) ~=~
  \sum_{z'\in\Z^d} \mu(z') 
\exp(a\cdot z') 
\]
is finite everywhere on $\R^d$ and hence, the condition (A2') is also satisfied. 
Moreover, let $\hat{k}$ be the 
greatest common divisor\footnote{Using the same arguments as in the proof of
  Lemma~\ref{lem1-4} one can easily show that $\hat{k}$ is the period of the
  random walk $(Z(t))$} of the set 
\[
{\cal K} ~=~ \left\{ n\in\N^* :~ \inf_{z\in\Z^{d-1}\times\N^*} \P_z(Z_+(n)=z)  ~>~
0\right\}. 
\]
If a non-negative function $\tilde{h}:\Z^{d-1}\times\N^*\to\R_+$ is harmonic for the Markov process
$(Z_+(t))$ and satisfies the equality  $\tilde{h}(z + \hat{k} w) = \tilde{h}(z)$ for all
$z\in\Z^{d-1}\times\N^*$ and $w\in\Z^{d-1}\times\{0\}$ then also 
\[
\sup_{x\in\Z^{d-1}} h(x,y) ~=~ \max_{x\in\Z^{d-1} :~ |x| \leq \hat{k} d} ~h(x,y)  ~<~ \infty 
\]
and consequently, by Proposition~\ref{pr5-1}, $\tilde{h} = c h_{0,+}$ for some
$c>0$. Hence, the condition (A3) is also satisfied with $h = h_{0,+}$. The hypotheses
 of Proposition~\ref{pr3-2} are therefore satisfied and  \eqref{e1-3} holds.  

To extend this result for an arbitrary $q\in{\cal S}_+^d$, it is sufficient to apply the above arguments
for a twisted random walk $\tilde{Z}(t) = (\tilde{X}(t),\tilde{Y}(t))$ on $\Z^d$ with transition probabilities 
\[
\tilde{p}(z,z') ~=~ \tilde\mu(z'-z) ~\dot=~ \mu(z'-z) \exp(a(q)\cdot(z'-z))
\]
for which  the mean is equal to 
\[
\tilde{m} ~=~ \sum_{z\in\Z^d} z \tilde\mu(z) ~=~ \nabla\varphi(a(q)).
\]
For such a random walk, $q=\tilde{m}/|\tilde{m}|\in{\cal S}^d_+$ and the hypotheses~(A) are satisfied if
they are satisfied for the initial random walk $(Z(t))$. Hence, using the same arguments as
above we conclude that Green's
function $\tilde{G}_+(z,z')$ of the random walk $(\tilde{Z}_+(t))$ killed outside of the
half-space $\Z^{d-1}\times\N^*$  satisfies the following property : for every
$z\in\Z^{d-1}\times\N^*$, 
\[
\lim_{n\to\infty} \tilde{G}_+(z,z_n)/\tilde{G}_+(z_0,z_n) ~=~
\tilde{h}_{0,+}(z)/\tilde{h}_{0,+}(z_0)
\]
where 
\[
\tilde{h}_{0,+}(z) ~\dot=~  \begin{cases}
y  -  \E_{z}\left( \tilde{Y}(\tau), \; \tau < +\infty
\right) &\text{if $\tilde{m}\in\R^{d-1}\times\{0\}$,}\\
\P_{z}\left( \tau = +\infty\right)
&\text{if $\tilde{m}\in\R^{d-1}\times]0,+\infty[$} 
\end{cases}
\]
with $\tau = \inf\{ t \geq 0 : \tilde{Y}(t) \leq 0\}$. 
Finally, using the identities 
\[
\tilde{G}_+(z,z') ~=~ e^{a(q)\cdot(z'-z)} G_+(z,z') \quad 
\text{
and } \quad 
\tilde{h}_{0,+}(z) ~=~ e^{-a(q)\cdot z} h_{a(q),+}(z) 
\]
we get \eqref {e1-3}. 
Theorem~\ref{th1-4} is therefore proved.

\section{Application for homogeneous random walk}\label{sec7}
In the present section, the Martin
boundary of a homogeneous random walk on $\Z^d$ is obtained as a consequence of
Proposition~\ref{pr3-2}. 

Remark that a random walk $Z(t)$ on $\Z^d$ with transition probabilities
$p(z,z')=\mu(z'-z)$, $z,z'\in\Z^d$, can be considered as a Markov additive process where the additive part
is $A(t)=Z(t)$ and the Markovian part is constant $M(t)\equiv 0$ with $E=\{0\}$.  In this
setting, Proposition~\ref{pr3-2} implies the following result.  
\medskip 

\noindent
{\bf Theorem~[Ney and Spitzer~\cite{Ney-Spitzer}]} {\em Suppose that the random walk 
  $Z(t)$  is irreducible with a non zero mean
\[ 
m ~=~ \sum_{z\in\Z^d} ~z~\mu(z) ~\not=~ 0, 
\]
and let the jump generating function $\varphi$ of \eqref{e1-0} 
be finite in a neighborhood of the set $D = \{a\in\R^d : \varphi(a)\leq 1\}$.  Then for
any $z\in\Z^d$, uniformly on $q\in{\cal S}^d$,
\[
G(z,z_n)/G(0,z_n) ~\to~ \exp(a(q)\cdot z) 
\]
when $|z_n|\to\infty$ and $z_n/|z_n|\to q$}
\begin{proof} Indeed, under the hypotheses of this theorem, the
  Markov additive process $(A(t),M(t))$ 
  with the additive part $A(t)=Z(t)$ and the constant Markovian part $M(t)\equiv 0$
  satisfies the conditions (A1') and (A2') (see the proof of Lemma~\ref{lem1-4} for more
  details). Hence, by Proposition~\ref{pr3-2}, for any sequence of points $z_n\in\Z^d$
  with $\lim_n|z_n| = \infty$, if
\be\label{e7p-4}
\liminf_{n\to\infty}\frac{1}{|z_n|} \log G(0, z_n) ~\geq~ 0 
\ee 
then also  
 \be\label{e7p-3}
\lim_{n\to\infty} G(z + \hat{k} w,z_n)/G(z,z_n) ~=~ 1, \quad \quad \forall  z,w\in\Z^d 
\ee
where $\hat{k}$ is the period of the random walk $(Z(t))$. Moreover,  the same arguments as in the proof
  of Proposition~\ref{pr4-2}\footnote{Here, instead of Proposition~\ref{pr4-1} one should use 
  Mogulskii's theorem} show that  
\[
\liminf_{n\to\infty}\frac{1}{|z_n|} \log G(0, z_n) ~\geq~ - ~a(q)\cdot q 
\]
when $|z_n|\to\infty$ and $z_n/|z_n| \to q\in{\cal S}^d$. Consider now a sequence of
  points $z_n\in\Z^d$ such that $|z_n|\to\infty$ 
  and $z_n/|z_n|\to q = q(0)\in{\cal S}^d$ as $n\to\infty$. Then $a(q) = 0$ and
  consequently, the
  last inequality provides  the inequality \eqref{e7p-4}. Hence, the equality  \eqref{e7p-3}
  holds. Moreover, let for a subsequence $(z_{n_k})$,  the sequence 
  of functions $G(\cdot, z_{n_k})/G(0, z_{n_k})$ converge pointwise in
  $\Z^d$ and let 
\[
h(z) ~=~ \lim_{k\to\infty} G(z, z_{n_k})/G(0, z_{n_k}). 
\]
Then  by Fatou lemma, the function $h$ is super-harmonic for $(Z(t))$ :
\[
\E_z\bigl(h(Z(t))\bigr) ~\leq~ h(z), \quad \quad \quad \forall \, t\in\N, \; z\in\Z^d
\]
and from \eqref{e7p-3} it follows that $
h(z+\hat{k}w) ~=~ h(z)$ for all  $z,w\in\Z^d$. 
Hence, $h$ is a super-harmonic  function with a finite number of
values.  By  Minimum principal (see~Woess\cite{Woess}) from this it follows that  the function $h$ is
constant. Moreover, since $h(0) = 1$ we conclude that $h(z)=1$ for all $z\in\Z^d$ and
consequently, 
\[
G(z, z_{n_k})/G(0, z_{n_k}) ~\to~ 1, \quad \quad \quad \forall
z\in\Z^d 
\]
when $|z_n|\to\infty$ and $z_n/|z_n| \to q(0)$ for any subsequence 
  $(z_{n_k})$ for which the sequence  
  of functions $G(\cdot, z_{n_k})/G(0, z_{n_k})$ converges pointwise in
  $\Z^d$. This proves that  
\[
G(z, z_{n})/G(0, z_{n}) ~\to~ 1, \quad \quad \quad \forall
z\in\Z^d 
\]
when $|z_n|\to\infty$ and $z_n/|z_n| \to q(0)$ because for every $z\in\Z^d$, the sequence $G(z,
z_{n})/G(0, z_{n})$ is bounded . 

Suppose now that $|z_n|\to\infty$
  and $z_n/|z_n|\to q \not= q(0)$ as $n\to\infty$. Then the same arguments applied for the
  twisted random walk having transition probabilities $
\tilde{p}(z,z') ~=~ \exp\bigl(a(q)\cdot(z'-z)\bigr) \mu(z'-z)$, the mean $
m(q) = \nabla\varphi(a(q))$ and Green's function 
$
\tilde{G}(z,z') = \exp\bigl(a(q)\cdot(z'-z)\bigr)G(z,z')$  
show that $\tilde{G}(z, z_{n})/\tilde{G}(0, z_{n}) ~\to~ 1$ as $n\to\infty$ for every $
z\in\Z^d$ 
and consequently 
\[
\lim_{n\to\infty}  G(z, z_{n})/G(0, z_{n}) ~=~ \exp(a(q)\cdot z), \quad \quad \quad \forall
z\in\Z^d.  
\]
Consider now the mapping $z\to w(z) = z/(1+|z|)$ from $\Z^d$ to ${\cal R} = \{w = z/(1+|z|) :
z\in\Z^d\}$ and its 
inverse mapping $w\to z(w) = w/(1-|w|)$.  Remark that the boundary $\partial{\cal
  R}$ of the set ${\cal R}$ in $\R^d$ is the $d$-dimensional 
sphere ${\cal S}^d$, the closure $\ol{\cal R}$ of the set ${\cal R}$ is compact and for
every $z\in\Z^d$, the function 
\[
f_z(w) ~=~ \begin{cases} G(z, z(w)))/G(0, z(w))) &\text{if $w\in{\cal R}$,}\\
\exp(a(w)\cdot z) &\text{if $w\in\partial{\cal R}$}
\end{cases}
\]
is continuous on $\ol{\cal R}$. Hence, for
every $z\in\Z^d$, the function $f_z$ is uniformly continuous on $\ol{\cal R}$ and
consequently, for any $\eps > 0$ there is $\delta > 0$ such that for all $q\in{\cal S}^d$
and $w\in{\cal R}$, 
\[
\left| G(z, z_n))/G(0, z_n))  - \exp(a(q)\cdot z) \right| ~<~ \eps 
\]
whenever $|q - w(z_n)| < \delta$.  This proves that $G(z,
z_{n})/G(0, z_{n}) ~\to~ \exp(a(q)\cdot z)$ when $|z_n|\to\infty$ and $z_n/|z_n| \to q$ 
uniformly on $q\in{\cal S}^d$ for every $z\in\Z^d$. 
\end{proof}

\section{Proof of Proposition~\ref{pr3-1}}\label{secA-1} 
To prove this proposition it is 
sufficient to show that 
\be\label{e8-1}
\limsup_{n\to\infty} ~{\cal G}(z+w',z_n)/{\cal G}(z+w,z_n) ~\leq~ 1
\ee
for any $z\in\Z^d\times E$. The main idea of the proof is
the following~:  Remark first of all that 
if \eqref{e3-4}  holds for some $z_0\in\Z^d\times E$ then also 
\be\label{e8-2}
\liminf_{n\to\infty}~\frac{1}{|z_n|} \log {\cal G}\bigl(z, z_n\bigr)  ~\geq~  0 
\ee
for all $z\in\Z^d\times E$ because the Markov-process ${\cal Z}(t)$ is irreducible. 
This inequality shows that the terms of the order $c\exp(-\delta |z_n|)$ give an
asymptotically negligible contribution to ${\cal G}(z,z_n)$. The quantities ${\cal
  G}(z,z_n)$ are decomposed into a negligible part
of the order $c\exp(-\delta |z_n|)$  and a main part for which it is possible to apply the method of
Bernoulli part decomposition developed in~\cite{McDonald:01,Foley-McDonald,McDonald:02}.

The first step of the proof shows that there is $\kappa >0$ for which the part 
\[
\sum_{0\leq t \leq \kappa |z_n|} \P_{z}\bigl({\cal Z}(t) = z_n\bigr) 
\]
of ${\cal G}(z,z_n)$ is negligible. Under the hypotheses of our lemma, there is $\delta_0
>0$ such that for any $0 < \delta \leq \delta_0$, 
\[
C_\delta ~\dot=~ \sup_{a\in\R^{d+k} : |a|\leq \delta} ~\sup_{z\in\Z^d\times E} ~\E_z(e^{a\cdot
  ({\cal Z}(1)-z)}) < \infty
\]
and hence, using Chebychev's inequality and Markov property we obtain that 
for any $l>0$, $z\in\Z^d\times E$ and any $a\in\R^d\times\R^k$ with $|a| \leq \delta \leq 
\delta_0 $, the following relation holds 
\begin{align*}
\P_{z}\left(a\cdot {\cal Z}(t) \geq  l\right) ~\leq~ e^{- l }
~\E_{z}\left( e^{a\cdot {\cal Z}(t)}  \right) ~\leq C_\delta^t   e^{- l + a\cdot z} ~\leq C_\delta^t   e^{- l +
  \delta |z|}, \quad \forall t\in\N.
\end{align*}
Using this inequality with $l= \delta |z_n| $ for $0 < \delta \leq
\delta_0$ we obtain 
\[
\P_{z}\bigl({\cal Z}(t) = z_n\bigr) ~\leq~ P_{z}\left(z_n\cdot {\cal Z}(t) = |z_n|^2\right) ~\leq~  C_\delta^t
\exp\left(- \delta |z_n| + \delta |z|\right) 
\]
and hence, for $\kappa = \delta/(2 \ln C_\delta)$ we get 
\begin{align}
\sum_{0\leq t \leq \kappa |z_n|} \P_{z}\bigl({\cal Z}(t) = z_n\bigr) 
&~\leq~ \exp\bigl(-\delta |z_n| + \delta |z|\bigl) \sum_{0\leq t \leq \kappa  |z_n| }
C_\delta^t \nonumber \nonumber  \\
&~\leq~ \exp\left(-  \frac{\delta}{2} |z_n| + \delta |z|\right)/(C_\delta-1) , \label{e8-3}
\end{align}
for any $z\in\Z^d\times E$ and any $n\in\N$. The left hand side of this inequality is
therefore a negligible part of ${\cal G}(z,z_n)$. 

Next, the method of Bernoulli part decomposition  is applied. 
Recall that under the hypotheses of our proposition, for  given $w,w'\in\Z^d\times\{0\}$, there
are $\eps > 0$ and $\hat{n} \in\N^*$ such that 
\be\label{e8-4}
\inf_{z\in\Z^d\times E} \min~\left\{ p^{(\hat{n})}(z,z+w), \, p^{(\hat{n})}(z,z+w')\right\} \geq \eps.
\ee
Suppose first that $\hat{n} = 1$ and let us consider two independent sequences $(\xi(t))$ and $(\zeta(t))$ of
independent  Bernoulli random 
variables with means $\E(\xi(t))=\eps$ and $\E(\zeta(t))=1/2$.   Then the Markov additive process
${\cal Z}(t) = (A(t),M(t))$ can be represented  in the 
following way~: at time $t$ we pick a  
random variable $\xi(t)$ and if $\xi(t)=1$, then we let ${\cal Z}(t+1) = {\cal Z}(t) +
\zeta(t) w + (1-\zeta(t)) w'$. Otherwise, for
each $y\in E$, we put $M(t+1) = y$ with probability  
\[
\Bigl(
p_M(M(t),y)\1_{\{M(t)\not=y\}} + (p_M(y,y)-\eps) \1_{\{M(t)=y\}}\Bigr)/(1-\eps) 
\]
where 
\[
p_M(y',y) \dot= \sum_{x\in\Z^d} p\bigl((0,y'),(x,y)\bigr), 
\]
and we let 
\begin{multline*}
A(t+1) ~=~ A(t) + \sum_{y\in E}  \1_{\{M(t)=M(t+1)=y\}}  B_{y,y}(t) \\ +
\sum_{y,y'\in E :y\not=y'}\1_{\{M(t)=y,M(t+1)=y'\}} B_{y,y'}(t) 
\end{multline*}
where $(B_{y,y'}(t), \; t\in\N, y\in E)$ is a family of mutually independent random
variables which are independent of the Markov process $(M(t))$ and of the sequences
$(\xi(t))$ and $(\zeta(t))$, such that for $x\in\Z^d$ and $y,y'\in E$,  
\[
\P(B_{y,y'}(t) = x) ~=~  {p\bigl((0,y),(x,y')\bigr)}/{p_M(y,y')} \quad
\text{when \; $y\not=y'$,}
\] 
\[
\P(B_{y,y}(t) = x) ~=~  \frac{p\bigl((0,y),(x,y)\bigr)}{(1-\eps)p_M(y,y)} \quad \text{if
  \; $(x,0)\not\in\{w,w'\}$,}
\]
and 
\[
\P(B_{y,y}(t) = x) ~=~  \frac{p\bigl((0,y),(x,y)\bigr) - \eps/2}{(1-\eps)p_M(y,y)} \quad \text{for
  \; $(x,0)\in\{w,w'\}$.}
\]
This representation shows that 
\[ 
{\cal Z}(t) ~=~ (A(t),M(t)) ~=~ \left(Q(t) + \sum_{s=1}^{N_t} \Bigl( \zeta(s) w + (1-\zeta(s)) w'\Bigr),
M(t)\right)  
\]
where $N_t = \text{Card}\{s\in\N : s\leq t, \, \xi(s) = 1\}$ is a Binomial random variable with mean
$\eps t$ and variance $\eps(1-\eps)t$, 
the random vector $(Q(t), M(t), N_t)$ is 
independent of the sequence $(\zeta(s))$ (the random variables  $Q(t)$, $M(t)$ and $N_t$ are
dependent) and the equality holds in a sense of the identity of
the distributions.  Hence, letting 
\[
L_n ~=~ \sum_{s=1}^n \zeta(s) 
\]
we obtain 
\begin{align}
\P_{z}&\left({\cal Z}(t) = z_n, \; N_t =N, \, L_N = l \; \right) \nonumber\\ 
&\quad  ~=~ \P_{z}\left( (Q(t),M(t)) = z_n - l w - (N-l)w', \; N_t = N, \; L_N = l \; \right) \nonumber\\
&\quad ~=~ \P_{z} \Bigl( (Q(t),M(t)) = z_n - l w - (N-l)w', \; N_t = N\Bigr) \, \P\left(
L_N = l \; \right) \label{e8-5}
\end{align} 
for all $z\in\Z^d\times E$, $n\geq 0$ and $0\leq l\leq N$. Furthermore, by Chebychev's inequality, 
\[
\P(N_t < \eps t/2) ~\leq~ \inf_{\eta < 0}  ~e^{- \eta \eps t/2} \E\left(e^{  \eta
  N_t}\right) ~=~ \exp\left( - t \theta_1\right)
\]
where $
\theta_1 ~\dot=~ \sup_{\eta < 0} \bigl( \eta \eps /2 - \log (\eps e^{\eta} +
  1-\eps)\bigr)  ~>~ 0$ 
because the function $f_1(\eta) = \eta \eps /2 - \log (\eps e^{\eta} +
  1-\eps)$ is concave, $f_1(0)=0$ and $f'_1(0) = -\eps/2 < 0$.  From this it follows that 
\begin{align}
\sum_{t >
  \kappa |z_n|}  \P_{z}({\cal Z}(t) = z_n, \; N_t < \eps t/2 ) 
&~\leq~ \sum_{t >
  \kappa |z_n|}  ~\P\left( N_t < \eps t/2\right) \nonumber \\ &~\leq~ \exp\left( - \kappa
  \theta_1 |z_n|\right)/(1-\exp( - \theta_1)) \label{e8-6} 
\end{align}
Hence, the left hand side of this  inequality is a negligible part of ${\cal
  G}(z,z_n)$. 
Moreover, for $0 < \sigma < 1/2$, 
\begin{align*}
\P\left(\left|L_N - \frac{N}{2 }\right| > \sigma N\right) &=~ 2 \P\left( L_N >  \frac{N}{2 } + \sigma N\right)\\
&\leq~ 2 ~\inf_{\eta > 0}  ~e^{- \eta (1/2 + \sigma )N} \E\left(\exp\left( \eta \sum_{s=1}^{N}
  \zeta(s)\right)\right)\\
&\leq 2 ~\exp\left(- N \theta_2\right)
\end{align*}
where $\theta_2 ~\dot=~ \sup_{\eta > 0}  \left(\eta (1/2 + \sigma ) - \log \bigl((e^\eta +
1)/2\bigr)\right) ~>~ 0$ 
because the function $f_2(\eta) = \eta (1/2 + \sigma ) - \log \bigl((e^\eta +
1)/2\bigr)$ is concave, $f_2(0)=0$ and $f_2'(0) = \sigma > 0$. Using this inequality we
get 
\begin{align}
\sum_{t >
  \kappa |z_n|} &\P_{z}\left({\cal Z}(t) = z_n, \; N_t \geq \eps t/2, \; \left|L_{N_t} - N_t/2\right| > \sigma N_t \; \right) \nonumber\\ 
&\leq \sum_{t >
  \kappa |z_n|}  ~\sum_{N\geq \eps t/2} ~\P\left( N_t = N, \,\left| L_N -
  N/2\right| > \sigma N\right) \nonumber\\ 
&\leq 2 \sum_{t >
  \kappa |z_n|}  ~\sum_{N \geq \eps t/2}~\exp\left(- N \theta_2\right) ~\leq~ \sum_{t >
  \kappa |z_n|}  \exp\left(-
  \eps \theta_2 t/2\right)/(1-\exp(- \theta_2))\nonumber\\
&\leq ~\exp\left(-
  \eps \theta_2 \kappa |z_n|/2\right)/ (1-\exp(- \theta_2))(1-\exp(-
  \eps \theta_2/2)). \label{e8-7} 
\end{align}
The last inequality shows that the part 
\[
\sum_{t >
  \kappa |z_n|} \P_{z}\left({\cal Z}(t) = z_n, \; N_t \geq \eps t/2, \; \left|L_{N_t} -
  N_t/2\right| > \sigma N_t \; \right)  
\]
of ${\cal G}(z,z_n)$ is also negligible. Finally, remark that for all those $l\in\N$ for which $|l-N/2| \leq \sigma
N$, the following relation holds  
\[
\P\left(L_N = l \right) ~=~  \P\left(L_N = l+1 \right) \frac{l+1}{N-l} ~\leq~ \P\left(L_N
= l+1 \right) \frac{1 + 2\sigma  +2/N}{1 - 2\sigma } 
\]
and hence, using  \eqref{e8-5} we obtain  
\begin{align*}
&\P_{z}\left({\cal Z}(t) = z_n, \; N_t =N, \, L_N = l \; \right) \\ 
&=~ \P_{z} \bigl( (Q(t),M(t)) = z_n - l w - (N-l)w', \; N_t = N\bigr) \, \P\left( L_N = l \; \right) \\
&\leq~ \frac{1 + 2\sigma  +2/N}{1 - 2\sigma } ~\P_{z}\bigl( (Q(t),M(t)) = z_n - l w -
  (N-l)w', \; N_t = N, \, L_N = l + 1\bigr) \\
&\leq~ \frac{1 + 2\sigma  +2/N}{1 - 2\sigma } ~\P_{z}\left({\cal Z}(t) = z_n + w - w', \; N_t =N,
  \, L_N = l + 1\; \right) 
\end{align*} 
From the last inequality it follows that   
\begin{align*}
\sum_{t >
  \kappa |z_n|} \P_{z}\bigl({\cal Z}(t) &= z_n, \; N_t \geq \eps t/2, \; \left| L_{N_t}-
  N_t/2\right| \leq \sigma N_t \; \bigr)  \nonumber\\ 
&~\leq~  \frac{1 + 2\sigma  +4/(\eps\kappa |z_n|)}{1 - 2\sigma } \sum_{t >
  \kappa |z_n|}  \P_{z}\left({\cal Z}(t) = z_n - w + w'\right) \nonumber \\ &~\leq~
  \frac{1 + 2\sigma  +4/(\eps\kappa |z_n|)}{1 - 2\sigma }  ~{\cal G}(z,z_n - w + w')
  \label{e8-8}   
\end{align*}
where ${\cal G}(z,z_n - w + w') = {\cal G}(z + w - w',z_n)$ because
$w,w'\in\Z^d\times\{0\}$. 
Using the above inequality together with \eqref{e8-3}, \eqref{e8-6} and \eqref{e8-7} 
we conclude
that for any $0 < \delta < \delta_0$ and $0 < \sigma < 1/2$ there are $\kappa > 0$, $C>0$ and
$\theta > 0$  such that      
\[
{\cal G}(z,z_n) ~\leq~ \frac{1 + 2\sigma
  +C/|z_n|}{1 - 2\sigma }~{\cal G}(z + w - w',z_n) + C \exp(- \theta |z_n| + \delta |z|). 
\]
for all $n\in\N$ and $z\in\Z^d\times E$. Because of \eqref{e8-2} from this it follows that 
\begin{align*}
\limsup_{n\to\infty} ~{\cal G}(z,z_n) / {\cal G}(z + w - w',z_n)
~\leq~ \frac{1 + 2\sigma}{1 - 2\sigma }, \quad \forall z\in\Z^d\times E 
\end{align*}
and consequently, also 
\begin{align*}
\limsup_{n\to\infty} ~{\cal G}(z+w',z_n) / {\cal G}(z + w,z_n)
~\leq~ \frac{1 + 2\sigma}{1 - 2\sigma }, \quad \forall z\in\Z^d\times E 
\end{align*}
(to get the last inequality it is sufficient to replace  $z$ by $z+w'$). Since $0 < \sigma <
1/2$ is arbitrary, then  letting at the last inequality $\sigma\to 0$ 
we get \eqref{e8-1}. When  \eqref{e8-4} holds with
$\hat{n}=1$, the inequality \eqref{e8-1}  is therefore verified.  

\medskip

Suppose now that $\hat{n} > 1$. Then for the Green's function
$\tilde{{\cal G}}(z,z_n) $ of the embedded Markov chain
$\tilde{\cal Z}(t) = (A(\hat{n}t),M(\hat{n}t))$ the above arguments
show that  for any $0 < \delta \leq \delta_0$ and $0 < \sigma < 1/2$ there are $\kappa >
0$, $C>0$, $\theta>0$  such that    
\[
\tilde{{\cal G}}(z,z_n)  ~\leq~ \frac{1 + 2\sigma
  +C/|z_n|}{1 - 2\sigma }~\tilde{{\cal G}}(z + w - w',z_n) + C \exp(- \theta |z_n| + \delta |z|). 
\] 
for all $z\in\Z^d\times E$. Using the 
identities 
\[
{\cal G}(z,z_n) ~=~ \sum_{t=0}^{\hat{n}-1} \sum_{z'} p^{(t)}(z,z')
\tilde{{\cal G}}(z',z_n) 
\]
and
\[
{\cal G}(z+ w-w',z_n) ~=~ \sum_{t=0}^{\hat{n}-1} \sum_{z'} p^{(t)}(z,z')
\tilde{{\cal G}}(z'+w-w',z_n) 
\]
it follows  therefore  that 
\begin{multline}\label{e8-8}
{\cal G}(z,z_n) \leq \frac{1 + 2\sigma  + C/|z_n|}{1 - 2\sigma }
{\cal G}(z+w-w',z_n)   + \\
C \sum_{t=0}^{\hat{n}-1} \sum_{z'} p^{(t)}(z,z') \exp\left(\delta |z'| - \theta |z_n|\right)
\end{multline}
Since under the hypotheses $(A_2)$, for $0 < \delta \leq \delta_0$,  
\begin{align*}
\sum_{z'} p^{(t)}(z,z') \exp\left(\delta |z'| \right) ~<~ \infty 
\end{align*}
then the inequality~\eqref{e8-8} combined with \eqref{e8-2} proves \eqref{e8-1}.  
Proposition~\ref{pr3-1} is  proved.

\bibliographystyle{amsplain}

\begin{thebibliography}{10}

\bibitem{Alili-Doney}
L.~Alili and R.~A. Doney, \emph{Martin boundaries associated with a killed
  random walk}, Ann. Inst. H. Poincar\'e Probab. Statist. \textbf{37} (2001),
  no.~3, 313--338.

\bibitem{Billingsley}
Patrick Billingsley, \emph{Convergence of probability measures}, Wiley series
  in probability and mathematical statistics, John Wiley \& Sons Ltd, New York,
  1968.

\bibitem{Borovkov:02}
A.A. Borovkov and A.A Mogulskii, \emph{The second rate function and the
  asymptotic problems of renewal and hitting the boundary for multidimensional
  random walks}, Siberian Math Journal \textbf{33} (1992), no.~4, 745--782.

\bibitem{Cartier}
P.~Cartier, \emph{Fonctions harmoniques sur un arbre}, Symposia Mathematica
  \textbf{9} (1972), 203--270.

\bibitem{McDonald:01}
B.~Davis and D.~McDonald, \emph{An elementary proof of the local central limit
  theorem.}, J.Theoret.Probab. \textbf{8} (1995), 693--701.

\bibitem{D-Z}
Amir Dembo and Ofer Zeitouni, \emph{Large deviations techniques and
  applications}, Springer-Verlag, New York, 1998.

\bibitem{Foley-McDonald}
Robert D.Foley and David R.McDonald, \emph{Bridges and networks: exact
  asymptotics}, Ann. Appl. Probab. \textbf{15} (2005), no.~1B, 542--586.

\bibitem{Doney:02}
R.~A. Doney, \emph{The {M}artin boundary and ratio limit theorems for killed
  random walks}, J.London Math.Soc. \textbf{2} (1998), no.~58, 761--768.

\bibitem{Doob}
J.~L. Doob, \emph{Discrete potential theory and boundaries}, J.Math. and Mech.
  \textbf{8} (1959), 433--458.

\bibitem{Dynkin:01}
E.B. Dynkin, \emph{Boundary theory of {M}arkov processes}, Russian Mathematical
  Surveys \textbf{24} (1969), no.~7, 1--42.

\bibitem{Hennequin}
P.L. Hennequin, \emph{Processus de {M}arkoff en cascade}, Ann. Inst. H.
  Poincar\'e \textbf{18} (1963), no.~2, 109--196.

\bibitem{Hunt}
G.~A. Hunt, \emph{Markoff chains and {M}artin boundaries}, Illinois J. Math.
  \textbf{4} (1960), 313--340.

\bibitem{Ignatiouk:02}
Irina Ignatiouk-Robert, \emph{Sample path large deviations and convergence
  parameters}, Annals of Applied Probability \textbf{11} (2001), no.~4,
  1292--1329.

\bibitem{Ignatiouk:08}
\bysame, \emph{Martin boundary of a reflected random walk on a half-space},
  2006, preprint: {\em http://arxiv.org/abs/math.PR/0610242}.

\bibitem{Kurkova-Malyshev}
I.A. Kurkova and V.A. Malyshev, \emph{Martin boundary and elliptic curves.},
  Markov Processes Related Fields \textbf{4} (1998), 203--272.

\bibitem{Martin}
R.S. Martin, \emph{Minimal positive harmonic functions}, Trans. Amer. Math.
  Soc. \textbf{49} (1941), 137--172.

\bibitem{McDonald:02}
D.~McDonald, \emph{On local limit theorems for integer valued random
  variables.}, Theory Probab. Appl. \textbf{24} (1979), 613--619.

\bibitem{Ney-Spitzer}
P.~Ney and Spitzer F, \emph{The {M}artin boundary for random walk}, Trans.
  Amer. Math. Soc. (1966), no.~121, 116--132.

\bibitem{R}
R.~Tyrrell Rockafellar, \emph{Convex analysis}, Princeton University Press,
  Princeton, NJ, 1997, Reprint of the 1970 original, Princeton Paperbacks.

\bibitem{Rogers:05}
L.~C.~G. Rogers and David Williams, \emph{Diffusions, {M}arkov processes, and
  martingales. {V}ol. 1: Foundations}, second ed., John Wiley \& Sons Ltd.,
  Chichester, 1994.

\bibitem{Rudin}
Walter Rudin, \emph{Real and complex analysis}, McGraw-Hill series in higher
  mathematics, McGraw-Hill book company, Ney-York, 1974.

\bibitem{Seneta}
E.~Seneta, \emph{Nonnegative matrices and {M}arkov chains}, second ed.,
  Springer-Verlag, New York, 1981.

\bibitem{Spitzer}
F.~Spitzer, \emph{Principles of random walk}, D. van Nostrand Company, INC,
  1964.

\bibitem{Woess}
Wolfgang Woess, \emph{Random walks on infinite graphs and groups}, Cambridge
  University Press, Cambridge, 2000.

\end{thebibliography}

\providecommand{\bysame}{\leavevmode\hbox to3em{\hrulefill}\thinspace}
\providecommand{\MR}{\relax\ifhmode\unskip\space\fi MR }
\providecommand{\MRhref}[2]{%
  \href{http://www.ams.org/mathscinet-getitem?mr=#1}{#2}
}
\providecommand{\href}[2]{#2}

\end{document}